\documentclass[11pt]{article}
\usepackage{amscd}
\usepackage{amsmath}
\usepackage{latexsym}
\usepackage{amsfonts}
\usepackage{amssymb}
\usepackage{amsthm}
\usepackage{graphicx}

 \newtheorem{proposition}{Proposition}[section]
 \newtheorem{definition}{Definition}[section]
 \newtheorem{lemma}{Lemma}[section]
 \newtheorem{theorem}{Theorem}[section]
 \newtheorem{corollary}{Corollary}[section]
 \newtheorem{remark}{Remark}[section]

\makeatletter
\newcommand{\Extend}[5]{\ext@arrow0099{\arrowfill@#1#2#3}{#4}{#5}}
\makeatother

\begin{document}

 \title{ On the blow up phenomenon for the mass critical \\ focusing Hartree
 equation in $ \mathbb{R}^4$
 }
\author{{Changxing Miao$^{\dag}$\ \ Guixiang Xu$^{\dag}$ \ \ and \ Lifeng Zhao $^{\ddag}$}\\
         {\small $^{\dag}$Institute of Applied Physics and Computational Mathematics}\\
         {\small P. O. Box 8009,\ Beijing,\ China,\ 100088}\\
         {\small $^\ddag$ Department of Mathematics, University of Science and Technology of China}\\
         {\small (miao\_changxing@iapcm.ac.cn, \ xu\_guixiang@iapcm.ac.cn, zhao\_lifeng@iapcm.ac.cn ) }\\
         \date{}
        }
\maketitle

\begin{abstract} We characterize the dynamics of the finite time blow up solutions with  minimal mass for the focusing
mass critical Hartree equation with $H^1(\mathbb{R}^4)$ data and
$L^2(\mathbb{R}^4)$ data, where we make use of the refined
Gagliardo-Nirenberg inequality of convolution type and the profile
decomposition.  Moreover, we also analyze the mass concentration
phenomenon of such blow up solutions.

\end{abstract}

\section{Introduction}
\setcounter{section}{1}\setcounter{equation}{0} In this paper, we
consider the Cauchy problem for the following Hartree equation
\begin{equation} \label{equ1}
\left\{ \aligned
    iu_t +  \Delta u  & = f(u), \quad  \text{in}\  \mathbb{R}^d \times \mathbb{R},\\
     u(0)&=u_0(x), \quad \text{in} \ \mathbb{R}^d.
\endaligned
\right.
\end{equation}
Here $f(u)=\lambda\big(V* |u|^2 \big)u$, $V(x)=|x|^{-\gamma},
0<\gamma<d$, and $*$ denotes the convolution in $\mathbb{R}^{d}$. If
$\lambda>0$, we call the equation $(\ref{equ1})$ defocusing; if
$\lambda<0$, we call it focusing. This equation describes the
mean-field limit of many-body quantum systems; see, e.g.,
\cite{FrL04}, \cite{Gi} and \cite{wiki}. An essential feature of
Hartree equation is that the convolution kernel $V(x)$ still retains
the fine structure of micro two-body interactions of the quantum
system. By contrast, NLS arise in further limiting regimes where
two-body interactions are modeled by a single real parameter in
terms of the scattering length. In particular, NLS cannot provide
effective models for quantum system with long-range interactions
such as the physically important case of the Coulomb potential
$V(x)\thicksim |x|^{-(d-2)}$ in $d\geq 3$, whose scattering length
is infinite.

There are many works on the global well-posedness and scattering of
equation $(\ref{equ1})$. For the defocusing case with
$2<\gamma<\min(4,d)$, J. Ginibre and G. Velo \cite{GiV} proved the
global well-posedness and scattering results in the energy space.
Later, K. Nakanishi \cite{Na} made use of a new Morawetz estimate to
obtain the similar results for the more general functions $V(x)$.
Recently, the authors proved the global wellposedness and scattering
for the defocusing, energy critical Hartree equation, see
\cite{MiXZ2} and \cite{MiXZ3}. The global wellposedness and
scattering of the focusing, energy critical Hartree equation can
refer to \cite{LiMZ08} and \cite{MiXZ5}. In this paper, we mainly
aim to characterize the dynamics of the finite time blow up
solutions with minimal mass for the focusing $L^2$-critical Hartree
equation with $H^1(\mathbb{R}^4)$ data and $L^2(\mathbb{R}^4)$ data.

Now we recall the related results about the focusing mass critical
Schr\"{o}dinger equation
\begin{equation}\label{equ12}
iu_t+\Delta u=-|u|^{\frac{4}{d}}u, \ \ u(0)=u_0,
\end{equation}
where $d$ is the spatial dimension. Equation (\ref{equ12}) is called
mass critical due to scaling invariance. If $u_0\in H^1$ is radial,
the mass concentration phenomena of the blow up solution was
observed near the blow-up time in \cite{MeT}. Later on, the radial
assumption was removed by M. Weinstein \cite{wein89} and Nawa
\cite{nawa}. For more detailed analysis of the blow up dynamic of
(\ref{equ12}), see \cite{Mer1}, \cite{Mer2}, \cite{MeR1},
\cite{MeR2}, \cite{MeR05} and the references therein. If $u_0$ only
lies in $L^2$, the situation seems quite different because we cannot
use the energy conservation law. The pioneering work in this
direction is due to J. Bourgain \cite{Bo} for $d=2$, where he proved
that there exists a blow-up time $T^*$,
$$\lim_{t\uparrow T^*}\sup_{\text{cubes}\ I\subset \mathbb{R}^2,\atop
\text{side}(I)<(T^*-t)^\frac{1}{2}}\Big(\int_{I}|u(t,x)|^2dx\Big)^\frac{1}{2}\geq
c(\|u_0\|_{L_x^2})>0,$$ where $c(\|u_0\|_{L_x^2})$ is a constant
depending on the mass of the initial data. A new proof can be found
in S. Keraani \cite{Ker2} by means of the profile decomposition in
\cite{MeV}. Bourgain's result was extended to dimension $d=1$ by R.
Carles and S. Keraani \cite{CaK} and to dimension $d\geq3$ by P.
B\'{e}gout and A. Vargas \cite{BeV}. Recently, R. Killip, T. Tao and
M. Visan \cite{tao2} established global well-posedness and
scattering for (\ref{equ12}) with radial data in dimension two and
mass strictly smaller then that of the ground state. Later R.
Killip, M. Visan and X. Zhang \cite{tao3} extended the results to
$d\geq3$. We dealt with the corresponding problem for the Hartree
equation in \cite{MiXZ4}.

This paper is devoted to the study of the blow up behavior of the
mass-critical Hartree equation in dimension four:
\begin{equation} \label{equ42}
\left\{ \aligned
    iu_t +  \Delta u  & = -(|x|^{-2}*|u|^2)u, \quad  \text{in}\  \mathbb{R}^4 \times \mathbb{R},\\
     u(0)&=u_0(x), \quad \text{in} \ \mathbb{R}^4.
\endaligned
\right.
\end{equation}
The corresponding free equation is
\begin{equation} \label{equ41}
\left\{ \aligned
    iu_t +  \Delta u  & = 0, \quad  \text{in}\  \mathbb{R}^4 \times \mathbb{R},\\
     u(0)&=u_0(x), \quad \text{in} \ \mathbb{R}^4.
\endaligned
\right.
\end{equation}
Note that $\gamma=2$ is the unique exponent which is mass-critical
in the sense that the natural scaling
\begin{equation*}
u_{\lambda}(t,x)=\lambda^{2}u(\lambda^2 t, \lambda x),
\end{equation*}
leaves the mass invariant. At the same time, $|x|^{-2}$ is just the
physically important case of Coulomb potential for dimension $d=4$.
Moreover, equation (\ref{equ42}) also possesses the pseudo-conformal
symmetry: If $u(t,x)$ solve (\ref{equ42}), then so does:
\begin{equation}\label{pcs}
\aligned v(t,x)=\frac{1}{|T-t|^2}\overline{u}(\frac{1}{t-T},
\frac{x}{t-T})e^{i\frac{|x|^2}{4(t-T)}}.
\endaligned
\end{equation}

 We firstly deal with equation (\ref{equ42}) with data in $H^1(\mathbb{R}^4)$. For the solution $u(t) \in H^1$ of (\ref{equ42}), there are the
following conserved quantities:
$$M(u(t))= \|u(t)\|_{L^2_x} =\|u(0)\|_{L^2_x} ,$$
$$E(u(t))= \frac{1}{2}\int_{\Bbb R^4}|\nabla
u|^2dx-\frac{1}{4}\int_{\Bbb R^4}\!\int_{\Bbb
R^4}\frac{|u(x)|^2|u(y)|^2}{|x-y|^2}dxdy = E(u(0)).$$
According to
the local wellposedness theory \cite{Ca03}, \cite{MiXZ1}, the
solution $u(t) \in H^1(\mathbb{R}^4)$ of (\ref{equ42}) blows up at
finite time $T$ if and only if
\begin{equation*}
\aligned \lim_{t\rightarrow T}\|\nabla
u(t)\|_{L^2}\rightarrow+\infty.
\endaligned
\end{equation*}

The blow-up theory is mainly connected to the notion of ground
state: the unique radial positive solution of the elliptic equation
\begin{equation}\label{ground}
 -\Delta Q + Q=(V*|Q|^2)Q.
\end{equation}
The existence of the positive solution is proved by the
concentration compactness principle at the beginning of Section 3,
which is close related to a refined Gagliardo-Nirenberg inequality
of convolution
type:
\begin{equation}
\label{equ3}\|u\|_{L^V}^4\leq\frac{2}{\|Q\|_{L^2}^2}\|u\|_{L^2}^2\|\nabla
u\|_{L^2}^2,
\end{equation}
where the definition of $L^V$ norm is given by (\ref{LV}). The
radial symmetry of the positive solution can be obtained from
\cite{Liu}. By adapting Lieb's uniqueness proof in \cite{lieb} for
the ground states $\phi\in H^1$ of the Choquard-Pekar equation
($V(x)=|x|^{-1}$ in  dimension $d=3$), the analogous result for
$(\ref{ground})$ can be obtained. See details in \cite{KrLR08}.
However, the uniqueness proof strongly depends on the specific
features of equation $(\ref{ground})$. It is different from the
corresponding results for semilinear elliptic equation in
\cite{kong}. As our result (Theorem \ref{thm4}) depends on the
uniqueness of the ground state of equation $(\ref{ground})$, it is
the reason why we do for the case $d=4$.

Together with the notion of the ground state $Q$, the invariance
$(\ref{pcs})$ yields an explicit blow-up solutions such that
$\big\|u\big\|_{L^2} =\big\|Q\big\|_{L^2}$. One can ask if there are
other finite time blow up solutions of $(\ref{equ42})$ with minimal
mass $\big\|Q\big\|_{L^2}$ and how to characterize the dynamics of
such blow up solutions near the blow up time.

Now, we can characterize the finite time blow-up solutions with
minimal mass in $H^1(\mathbb{R}^4)$.
\begin{theorem}\label{thm4}Let $u_0\in H^1(\Bbb R^4)$ such that $\|u_0\|_{L^2}=\|Q\|_{L^2}$ and $u$ be the blow up solution of (\ref{equ42})
at finite time $T$,  then there exists $x_0\in\mathbb{R}^4$ such
that $e^{i\frac{|x-x_0|^2}{4T}}u_0\in\mathcal{A}$, where
\begin{equation*}
\aligned
 \mathcal{A}=\bigg\{\rho^{2}e^{i\theta}Q(\rho x+y),
y\in\mathbb{R}^4, \rho\in\mathbb{R}_*^+, \theta\in[0,2\pi)\bigg\}.
\endaligned
\end{equation*}
\end{theorem}
\begin{theorem}\label{thm5}Let $u$ be a solution of (\ref{equ42})
which blows up at finite time $T>0$ with initial data $u_0\in
H^1(\Bbb R^4)$, and $\lambda(t)>0$ such that $\lambda(t)\|\nabla
u\|_{L^2}\rightarrow+\infty$ as $t\uparrow T$. Then there exists
$x(t)\in\mathbb{R}^4$ such that
$$\liminf_{t\uparrow T}\int_{|x-x(t)|\leq\lambda(t)}|u(t,x)|^2dx\geq\int_{\Bbb R^4}|Q|^2dx.$$\end{theorem}

The corresponding result of Theorem \ref{thm4} for the Schr\"odinger
equation has been established by F. Merle in \cite{Mer2}. The
corresponding result for Theorem \ref{thm5} was proved by M.
Weinstein in \cite{wein89}. T. Hmidi and S. Keraani gave a direct
and simplified proof of the above results in \cite{HMKer}. The new
ingredient for the Hartree equation is the refined
Gagliardo-Nirenberg inequality of the convolution type (\ref{equ3}),
whose proof is based on the well-known concentration compactness
method and thus one has to deal with the intertwining of convolution
and orthogonality.

Next we consider the blow up behavior of (\ref{equ42}) with $L^2$
data. In \cite{MiXZ1}, we showed that for any $u_0\in L^2(\Bbb
R^4)$, there exists a unique maximal solution $u$ to (\ref{equ42}),
with
$$u\in C((-T_*, T^*), L^2(\Bbb R^4))\cap L_{loc}^{3}((-T_*,
T^*), L^3(\Bbb R^4)), \ \
$$
and we have the following alternative: either $T_*=T^*=+\infty$ or
$$\min\{T_*,T^*\}<+\infty\ \ \text{and}\ \
\big\|u\big\|_{L_{t}^{3}((-T_*, T^*),\ L^3_x)}=+\infty.$$ Moreover,
there exists $\delta>0$ such that if
\begin{equation}\label{ap}\|u_0\|_{L^2}<\delta,\end{equation}
 the initial value problem (\ref{equ42}) has a unique global
 solution $u(t,x)\in L_{t,x}^3(\mathbb{R}\times\mathbb{R}^4)$. We define $\delta_0$ as the supremum of $\delta$ in (\ref{ap}) such
that the global existence for Cauchy problem (\ref{equ42}) holds,
with $u\in (C\cap L^\infty)(\Bbb R, L^2(\Bbb R^4))\cap L^3(\Bbb
R\times\Bbb R^4)$. Then in the ball $B_{\delta_0}:=\{u_0,
\|u_0\|_{L^2}<\delta_0\}$,
  (\ref{equ42}) admits a complete scattering theory with respect
  to the associated linear problem.  Similar to the focusing mass-critical Schr\"{o}dinger equation,
  we also conjecture that $\delta_0$ should be $\|Q\|_{L^2}$ for the Hartree equation. We have verified the
  conjecture for radial data in \cite{MiXZ4}. For general data, it remains open.
  \begin{definition}Let $u_0\in L^2(\mathbb{R}^4)$. A solution of (\ref{equ42}) is said to be a blow-up solution for
  $t>0$, if $\ T^*<+\infty$ or
  \begin{equation*}
  T^*=+\infty\ \ \text{and}\ \ \|u\|_{L_{t}^3((0,
+\infty),\ L^3_x)}=+\infty.\end{equation*}
   Similarly for $t<0$.\end{definition}
   Now we are in position to state the existence of the blow up solutions in both time directions with minimal mass in $L^2(\mathbb{R}^4)$.
\begin{theorem}\label{thm-mass1}
There exists an initial data $u_0\in L^2(\mathbb{R}^4)$ with
$\|u_0\|_{L^2}=\delta_0$, for which the solution of (\ref{equ42})
blows up for both $t>0$ and $t<0$.
\end{theorem}

As a direct consequence of the above theorem and the
pseudo-conformal transform $(\ref{pcs})$, we obtain the existence of
the finite time blow up solutions with minimal mass in
$L^2(\mathbb{R}^4)$.

\begin{corollary}
There exists an initial data $u_0\in L^2(\mathbb{R}^4)$ with
$\|u_0\|_{L^2}=\delta_0$, for which the solutions of (\ref{equ42})
blows up at finite time $T^*>0$.
\end{corollary}

\begin{theorem}\label{thm-mass2}Let $u$ be a blow up solution
of  (\ref{equ42}) at finite time $T^*>0$ such that
$\|u_0\|_{L^2}<\sqrt{2}\delta_0$. Let $\{t_n\}_{n=1}^{\infty}$ be
any time sequence such that $t_n\uparrow T^*$ as
$n\rightarrow\infty$, and let $\lambda(t)>0$, such that
$$\frac{\sqrt{T^*-t}}{\lambda(t)}\rightarrow0,\ \text{as}\ t\uparrow T^*.$$ Then there exist a subsequence of
$\{t_n\}_{n=1}^{\infty}$ {\rm(}still denoted by $\{t_n\}${\rm)} and
$x(t)\in \mathbb{R}^4$ that satisfy the following properties.
\begin{enumerate}
\item[$(i)$] There exists a function $\psi\in L^2(\Bbb R^4)$ with
$\|\psi\|_{L^2}\geq\delta_0$ such that the solution $U$ of
(\ref{equ42}) with initial data $\psi$ blows up for both $t>0$ and
$t<0$.

\item[$(ii)$] There exists a sequence $\{\rho_n, \xi_n,
x_n\}_{n=1}^{\infty}\subset \mathbb{R}^*_+\times \mathbb{R}^4\times
\mathbb{R}^4$ such that
\begin{equation*}
\rho_n^{2}e^{ix\cdot\xi_n}u(t_n,\rho_nx+x_n)\rightharpoonup\psi,\ \
\text{weakly in}\ \ L^2.
\end{equation*}
Furthermore, we have
\begin{equation*}\lim_{n\rightarrow\infty}\frac{\rho_n}{\sqrt{T^*-t_n}}\leq\frac{1}{\sqrt{T^{**}}}
\end{equation*}
where $T^{**}$ denotes the lifespan of $U$.

\item[$(iii)$] $$\liminf_{t\uparrow
T^*}\int_{|x-x(t)|\leq\lambda(t)}|u(x,t)|^2dx\geq\delta_0^2.$$
\end{enumerate}
\end{theorem}

\begin{corollary}\label{co} Let $u$ be a blow up solution with minimal mass
of (\ref{equ42}) at finite time $T^*>0$. Let $\{t_n\}_{n=1}^\infty$
be any time sequence such that $t_n\uparrow T^*$ as
$n\rightarrow\infty$. Then there exists a subsequence of
$\{t_n\}_{n=1}^\infty$ {\rm(} still denoted by
$\{t_n\}_{n=1}^\infty$ {\rm)} and $x(t)\in \mathbb{R}^4$ that
satisfy the following properties:
\begin{enumerate}
\item[$(i)$] There exists a function $\psi\in L^2(\Bbb R^4)$ with
$\|\psi\|_{L^2}\geq\delta_0$ such that the solution $U$ of
(\ref{equ42}) with initial data $\psi$ blows up for both $t>0$ and
$t<0$.

\item[$(ii)$] There exists a sequence $\{\rho_n, \xi_n,
x_n\}_{n=1}^{\infty}\subset \mathbb{R}^*_+\times \mathbb{R}^4\times
\mathbb{R}^4$ such that
\begin{equation*}
\rho_n^{2}e^{ix\cdot\xi_n}u(t_n,\rho_nx+x_n)\rightarrow \psi,\ \
\text{strongly in}\ \ L^2.
\end{equation*}
Furthermore, we have
\begin{equation*}\lim_{n\rightarrow\infty}\frac{\rho_n}{\sqrt{T^*-t_n}}\leq\frac{1}{\sqrt{T^{**}}}
\end{equation*}
where $T^{**}$ denotes the lifespan of $U$.

\item[$(iii)$] $$\liminf_{t\uparrow
T^*}\int_{|x-x(t)|\leq\lambda(t)}|u(x,t)|^2dx\geq\delta_0^2.$$
\end{enumerate}
\end{corollary}

Similar results for the nonlinear Schr\"odinger equation have
appeared in F. Merle, L. Vega \cite{MeV} and S. Keraani \cite{Ker2}.
Since the nonlinearity is non-local for the Hartree equation, we
have to pursue suitable decomposition in physical space to exploit
the orthogonality.

 We will often use the notations $a\lesssim b$ and $a=O(b)$ to
mean that there exists some constant $C$ such that $a\leq Cb$. The
derivative operator $\nabla$ refers to the derivative with respect
to space variable only. We also occasionally use subscripts to
denote the spatial derivatives and use the summation convention over
repeated indices.

For $1\leq p\leq \infty,$ we define the dual exponent $p'$ by
$\frac1p+ \frac1{p'} =1$. For any time interval $I$, we use
$L^q_tL^r_x(I\times \mathbb{R}^4)$ to denote the spacetime Lebesgue
norm
\begin{equation*}
\big\|u\big\|_{L^q_tL^r_x(I\times \mathbb{R}^4)}:=\bigg(
\int_I\big\| u\big\|^q_{L^r(\mathbb{R}^4)}dt \bigg)^{1/q}
\end{equation*}
with the usual modifications when $q=\infty$. When $q=r$, we
abbreviate $L^q_tL^r_x$ by $L^q_{t,x}$.

We say that a pair $(q, r)$ is admissible if
\begin{equation*}
\frac{2}{q} = 4\Big(\frac{1}{2}-\frac{1}{r}\Big),\ \ 2\leq q\leq
+\infty.
\end{equation*}
For a spacetime slab $I\times \mathbb{R}^4$, we define the {\it
Strichartz} norm $\dot{S}^0(I)$ by
\begin{equation*}
\big\|u\big\|_{\dot{S}^0(I)}:= \sup_{(q, r)\ \text{admissible}}
\big\|u\big\|_{L^q_tL^r_x(I\times \mathbb{R}^4)}.
\end{equation*}
and define $\dot{S}^1(I)$ by
\begin{equation*}
\big\|u\big\|_{\dot{S}^1(I)}:= \big\|\nabla
u\big\|_{\dot{S}^0(I)}.
\end{equation*}
We also define $\dot{\mathcal{N}}^0$ as the Banach dual space  of
$\dot{S}^0$.

 Throughout this paper, we denote
\begin{equation}\label{LV}\big\|u\big\|_{L^V}:=\Big(\int\!\!\int
|u(x)|^{2}V(x-y)|u(y)|^{2}dxdy\Big)^{\frac{1}{4}}.\end{equation}

The rest of this paper is organized as follows: In Section 2, we
 recall the preliminary estimates such as Strichartz estimates and
 Virial identity. In Section 3, we prove Theorem \ref{thm4}
 and Theorem \ref{thm5}. Section 4 is devoted to the proof of  Theorem \ref{thm-mass1} and Theorem
 \ref{thm-mass2}.

\section{Preliminaries}
\setcounter{equation}{0} We now recall some useful estimates. First,
we have the following {\it Strichartz } inequalities

\begin{lemma}[\cite{Ca03}, \cite{KeT98}]\label{lemstri1}
Let $u$ be an $\dot{S}^0(I)$ solution to the Schr\"{o}dinger
equation in $(\ref{equ1})$. Then
\begin{equation*}
\big\|u\big\|_{\dot{S}^0} \lesssim \big\|u(t_0)\big\|_{L^2(\Bbb
R^4)} + \big\|f(u)\big\|_{L^{q'}_tL^{r'}_x(I\times \mathbb{R}^4)}
\end{equation*}
for any $t_0 \in I$ and any admissible pairs $(q, r)$. The
implicit constant is independent of the choice of interval $I$.
\end{lemma}

By definition, it immediately follows that for any function $u$ on
$I\times \mathbb{R}^4$,
\begin{equation*}
\big\|u\big\|_{L^{\infty}_tL^2_x} + \big\| u\big\|_{L_{t,x}^3}
\lesssim \big\| u\big\|_{\dot{S}^0},
\end{equation*}
where all spacetime norms are taken on $I\times \mathbb{R}^4$.

\begin{lemma}\label{nle} Let $\displaystyle
f(u)(t,x)=\pm u\big(V*|u|^2\big)(t,x)$, where $V(x)=|x|^{-2}$. For
any time interval $I$ and $t_0 \in I$, we have
\begin{equation*}
\Big\|\int^t_{t_0} e^{i(t-s)\Delta}f(u)(s,x) ds\Big\|_{\dot{S}^0(I)}
\lesssim \big\|u\big\|^3_{L_{t,x}^3}.
\end{equation*}
\end{lemma}
{\noindent \it Proof. } By Strichartz estimate,
Hardy-Littlewood-Sobolev inequality and H\"{o}lder inequality, we
have
\begin{equation*}
\aligned \Big\|\int^t_{t_0} e^{i(t-s)\Delta}f(u)(s,x)
ds\Big\|_{\dot{S}^0(I)} & \lesssim \|f(u)(t,x)\|_{L^1_tL^2_x}\\
& \lesssim
\|V*|u|^2\|_{L^\frac{3}{2}_tL^6_x}\|u\|_{L_{t,x}^3}\\
&\lesssim\|u\|_{L_{t,x}^3}^3.
\endaligned
\end{equation*}

In addition, we have obtained the Virial identity in the proof of
the localized Morawetz estimates \cite{MiXZ2}. Indeed, let
$V^a_0(t)=\displaystyle \int a(x)|u(t,x)|^{2}dx$, where $a(x)$ is
real-valued and $u$ is the solution of (\ref{equ1}) with
$f(u)=-\big(|x|^{-\gamma} \ast |u|^2\big) u$. Then we get
$$M_{0}^{a}(t)=:\partial_{t}V^a_0(t)=2\Im\int a_{j}u_{j}\overline{u}dx$$
and
\begin{equation}\label{symvirial}
\aligned
\partial_{t}M_{0}^{a}(t)=&-2\Im\int
a_{jj}u_{t}\overline{u}dx-4\Im\int a_{j}\overline{u}_{j}u_{t}dx\\
=&-\int\triangle\triangle a|u|^{2}dx+4\Re\int
a_{jk}\overline{u_{j}}u_{k}dx \\
& -\iint \big(\nabla a(x)- \nabla a(y)\big)\nabla
V(x-y)|u(y)|^2|u(x)|^{2}dxdy.\\
\endaligned
\end{equation}

\begin{lemma} If we choose $a(x)=|x|^2$, then we have
\begin{equation}\label{est1}\partial_{t}M^{a}_0(t)=8\int|\nabla u|^2dx-2\gamma\iint
V(x-y)|u(y)|^2|u(x)|^{2}dxdy.\end{equation}
\end{lemma}

\begin{lemma}
If $a(x)=|x|^2$ and $\gamma=2$, we have
\begin{equation}\label{equ21}\partial_{t}^2V_0^a(t)=16E(u(0)).\end{equation}
If $E(u(0))<0$,
 the nonnegative function $V_0^a(t)$ is concave, so the maximal interval of existence is finite.
 This yields that the solution of (\ref{equ42}) must blow up in both directions.
 \end{lemma}

  \section{The blow-up dynamics of the focusing mass critical Hartree equation with $H^1$ data}
  \setcounter{equation}{0}  Let $V(x)=|x|^{-2}$, we study the minimizing functional
\begin{equation*}
J:=\min\{J(u): u\in
  H^1(\mathbb{R}^4)\}, \quad \text{where} \quad J(u):=\frac{\|u\|_{L^2}^2\|\nabla
  u\|_{L^2}^2}{\|u\|_{L^V}^4}.
\end{equation*}

First, we have
\begin{lemma}\label{lem42}
If $W$ is the minimizer of $J(u)$, then $W$ satisfies
\begin{equation}\label{equ43}
\Delta W+\alpha(|x|^{-2}*|W|^2)W=\beta W,\quad \text{where}\quad
\alpha=\frac{2J}{\|W\|_{L^2}^2};\ \ \beta=\frac{\|\nabla
  W\|_{L^2}^2}{\|W\|_{L^2}^2}.
\end{equation}
\end{lemma}
\begin{remark}\label{rem}
If $W$ is minimizer of $J(u)$, then $|W|$ is also a minimizer.
Hence, we can assume that $W$ is positive. In fact, we have
$$-|\nabla W|\leq\nabla |W|\leq|\nabla W|$$ in the sense of
distribution. In particular, $|W|\in H^1$ and $J(|W|)\leq J(W)$.
\end{remark}

   {\noindent \it Proof of Lemma \ref{lem42}. } It follows from the fact that $W$, the minimizing function, is in $H^1(\mathbb{R}^4)$ and satisfies the
  Euler-Lagrange equation:
\begin{equation*}
  \frac{d}{d\varepsilon}J(W+\varepsilon
  v)\Big|_{\varepsilon=0}=0.
  \end{equation*}
  Equivalently, we have
  \begin{align*}&\|\nabla W\|_{L^2}^2\|W\|_{L^V}^4\int2\Re(W\bar{v})dx+\|W\|_{L^2}^2\|W\|_{L^V}^4\int2\Re(\nabla W\nabla
  \bar{v})dx\\
  & \ \ -\|\nabla
  W\|_{L^2}^2\|W\|_{L^2}^2\Big(\int(V*2\Re(W\bar{v}))|W|^2dx+\int(V*|W|^2)2\Re(W\bar{v})dx\Big)=0.
  \end{align*}

  Since
  $$\int(V*2\Re(W\bar{v}))|W|^2dx=\int(V*|W|^2)2\Re(W\bar{v})dx,$$
we have $$\Delta
W+\frac{2J}{\|W\|_{L^2}^2}(V*|W|^2)W=\frac{\|\nabla
W\|_{L^2}^2}{\|W\|_{L^2}^2}W.$$

\begin{proposition}\label{pro41}
J is attained at a function u with the following properties:
$$u(x)=aQ(\lambda x+b),\ \text{for some}\  a\in\mathbb{C}^*,\ \lambda>0, \text{and any}\  b\in\mathbb{R}^4.$$
where $Q$ satisfies (\ref{ground}). Moreover,
$$J=\frac{\|Q\|_{L^2}^2}{2}.$$
\end{proposition}

We prove this proposition by the following profile decomposition.
\begin{lemma}[Profile decomposition \cite{HMKer}]\label{lem43}
For a bounded sequence $\{u_n\}_{n=1}^{\infty}\subset
H^1(\mathbb{R}^4)$, there is a subsequence of
$\{u_n\}_{n=1}^{\infty}$ {\rm(}still denoted by $\{u_n\}${\rm)} and
a sequence $\{U^{(j)}\}_{j\geq1}$ in $H^1(\mathbb{R}^4)$ and for any
$j\geq1$, a family $(x_n^j)$ such that
\begin{enumerate}
\item[$(i)$] If $j\neq k$, $|x_n^j-x_n^k|\rightarrow\infty$, as
$n\rightarrow\infty$.
\item[$(ii)$] For every $l\geq1$,
\begin{equation}\label{equ44}u_n(x)=\sum_{j=1}^lU^{(j)}(x-x_n^j)+r_n^l(x).\end{equation} Moreover,
for any $p\in(2,4)$,
\begin{equation}\label{equ45}
\limsup_{n\rightarrow\infty}\|r_n^l\|_{L^p(\mathbb{R}^4)}\rightarrow0\
\ \text{as}\ \  l\rightarrow+\infty.
\end{equation}
\item[$(iii)$]\begin{equation}\label{equ46}\|u_n\|_{L^2}^2=\sum_{j=1}^l\|U^{(j)}\|_{L^2}^2+\|r_n^l\|_{L^2}^2+o_n(1),\end{equation}
\begin{equation}\label{equ47}\|\nabla u_n\|_{L^2}^2=\sum_{j=1}^l\|\nabla U^{(j)}\|_{L^2}^2+\|\nabla r_n^l\|_{L^2}^2+o_n(1).\end{equation}
\end{enumerate}
\end{lemma}

{\noindent \it Proof of Proposition \ref{pro41}.}  Choose  a
sequence $\{u_n\}_{n=1}^{\infty}\subset H^1(\mathbb{R}^4)$ such that
$J(u_n)\rightarrow J$. Suppose $\|u_n\|_{L^2}=1$ and
$\|u_n\|_{L^V}=1$, then

 $$J(u_n)=\int|\nabla u_n|^2dx\rightarrow
J.$$
Note that $\{u_n\}_{n=1}^{\infty}$ is bounded in $H^1$, by
Lemma \ref{lem43}, we have (\ref{equ44})-(\ref{equ47}). From
(\ref{equ46}) and (\ref{equ47}), we have
\begin{equation}\label{equ13}\sum_{j=1}^l\|U^{(j)}\|_{L^2}^2\leq1,\quad \
\sum_{j=1}^l\|\nabla U^{(j)}\|_{L^2}^2\leq J.\end{equation}
Moreover, by H\"older and Young  inequalities, we have
$$\|r_n^l\|_{L^V}^4\leq\|r_n^l\|_{L^\frac83}^4.$$ From
(\ref{equ45}),
$\displaystyle\limsup_{n\rightarrow\infty}\|r_n^l\|_{L^\frac83}\buildrel{l\rightarrow\infty}\over{\longrightarrow}0$.
It follows that
$$\displaystyle\limsup_{n\rightarrow\infty}\|r_n^l\|_{L^V}\buildrel{l\rightarrow\infty}\over{\longrightarrow}0.$$
Moreover,
\begin{align}
&\iint \frac{|\sum_{j=1}^{l}U^{(j)}(x-x_n^j)|^2|\sum_{j=1}^{l}U^{(j)}(y-x_n^j)|^2}{|x-y|^2}dxdy\nonumber\\
\leq&\sum_{j=1}^{l}\iint \frac{|U^{(j)}(x-x_n^j)|^2|U^{(j)}(y-x_n^j)|^2}{|x-y|^2}dxdy\label{equ48}\\
&+\sum_{j=1}^l\sum_{k\neq j}\iint
\frac{|U^{(j)}(x-x_n^j)||U^{(k)}(x-x_n^k)|
(\sum_{i=1}^{l}|U^{(i)}(y-x_n^i)|)^2}{|x-y|^2}dxdy\label{equ49}\\
&+\sum_{j=1}^l\sum_{k\ne j}\iint
\frac{|U^{(j)}(y-x_n^j)||U^{(k)}(y-x_n^k)|
(\sum_{i=1}^{l}|U^{(i)}(x-x_n^i)|)^2}{|x-y|^2}dxdy\label{equ499}\\
&+\sum_{j=1}^l\sum_{k\neq j}\iint
\frac{|U^{(j)}(x-x_n^j)|^2|U^{(k)}(y-x_n^k)|^2}{|x-y|^2}dxdy.\label{equ410}
\end{align}
Without loss of generality we can assume that all $U^{(j)}$'s are
continuous and compactly supported. Then
$$(\ref{equ48})=\sum_{j=1}^{l}\iint \frac{|U^{(j)}(x)|^2|U^{(j)}(y)|^2}{|x-y|^2}dxdy,$$
and by orthogonality, we have \begin{eqnarray*} \aligned
(\ref{equ49})\leq\sum_{i=1}^{l} \sum_{j=1}^l\sum_{k\neq j}
\|U^{(i)}(y-x_n^i)\|_{L^{\frac83}}^2
\|U^{(j)}(\cdot-x_n^j)U^{(k)}(\cdot-x_n^k)\|_{L^\frac{4}{3}}\rightarrow0,
\quad n\longrightarrow\infty.
\endaligned
\end{eqnarray*}
(\ref{equ499}) can be similarly estimated. At last,  we estimate
\begin{eqnarray*}
\aligned
(\ref{equ410})&=\sum_{j=1}^{l}\sum_{k\neq j}\iint \frac{|U^{(j)}(x)|^2|U^{(k)}(y)|^2}{|x-y-x_n^j+x_n^k|^2}dxdy\\
&\leq \sum_{j=1}^l\sum_{k\neq j}
\frac{C}{|x_n^j-x_n^k|^2}\|U^{(j)}\|_{L^2}^2\|U^{(k)}\|_{L^2}^2
\rightarrow 0,  \quad n\longrightarrow\infty.
\endaligned
\end{eqnarray*}
Therefore, we conclude
$$\Big\|\sum_{j=1}^lU^{(j)}(x-x_n^j)\Big\|_{L^V}^4\rightarrow\sum_{j=1}^l\|U^{(j)}\|_{L^V}^4\quad \text{as}\quad n\longrightarrow\infty.$$
Thus, we have
\begin{equation*}
\lim_{l\rightarrow\infty}\sum_{j=1}^l\|U^{(j)}\|_{L^V}^4=1.
\end{equation*}

By the definition of $J$, we have
$$J\|U^j\|_{L^V}^4\leq\|U^{(j)}\|_{L^2}^2\|\nabla U^{(j)}\|_{L^2}^2.$$
So we get that
\begin{equation*}
J\sum_{j=1}^l\|U^j\|_{L^V}^4\leq\sum_{j=1}^l\|U^{(j)}\|_{L^2}^2\|\nabla
U^{(j)}\|_{L^2}^2.
\end{equation*}

On the other hand,
$$\sum_{j=1}^l\|U^{(j)}\|_{L^2}^2\|\nabla U^{(j)}\|_{L^2}^2\leq
\sum_{j=1}^l\|U^{(j)}\|_{L^2}^2\sum_{j=1}^l\|\nabla
U^{(j)}\|_{L^2}^2\leq J.$$ Thus we conclude that only one term
$U^{(j_0)}$ is non-zero, i. e.
\begin{equation}\label{equ421}\|U^{(j_0)}\|_{L^2}=1;\ \ \|U^{(j_0)}\|_{L^V}=1;\ \
\|\nabla U^{(j_0)}\|_{L^2}^2=J.\end{equation} This shows that
$U^{(j_0)}$ is the minimizer of $J(u)$. From (\ref{equ421}), we have
$$\Delta U^{(j_0)}+2J(|x|^{-2}*|U^{(j_0)}|^2)U^{(j_0)}=JU^{(j_0)}.$$
By Remark \ref{rem}, we can assume that $U^{j_0}$ is positive. Let
$U^{(j_0)}=aQ(\lambda x+b)$, where $Q$ is the positive solution of
(\ref{ground}). An easy computation gives that $\lambda^2=2a^2=J$.

 Next we compute the best constant
$J$ in terms of $Q$. Multiplying  $(\ref{ground})$ by $Q$ and
integrating both sides of
 this equation, we have
\begin{equation}\label{equ52}-\int|\nabla
Q|^2dx+\int(V*|Q|^2)|Q|^2dx=\int|Q|^2dx.\end{equation} Since
$$\int (x\cdot\nabla Q)Qdx=-2\int |Q|^2dx,$$
\begin{align*}
\int x\cdot\nabla Q\Delta Qdx&=-\sum_{i,j}\int \big(
\delta_{ij}\partial_iQ\partial_jQ+
x_i\partial_i\partial_jQ\partial_jQ\big)=\|\nabla Q\|_{L^2}^2,
\end{align*}
and \begin{align*}\int x\cdot\nabla Q(V*|Q|^2)Qdx&=\frac{1}{2}\int
x\cdot\nabla Q^2(V*|Q|^2)dx\\
&=\frac{1}{2}\int
x\cdot\nabla\big((V*|Q|^2)Q^2\big)dx-\frac{1}{2}\int x\cdot(\nabla
V*Q^2)Q^2dx\\
&=-2\int(V*|Q|^2)Q^2dx+\iint \frac{x\cdot(x-y)}{|x-y|^4}Q^2(x)Q^2(y)dxdy\\
&=-\frac32\|Q\|_{L^V}^4,\end{align*}
 we have
$$ \|\nabla Q\|_{L^2}^2-\frac32\|Q\|_{L^V}^4=-2\|Q\|_{L^2}^2.$$
 Together with (\ref{equ52}), this yields
$$\|\nabla Q\|_{L^2}^2=\|Q\|^2_{L^2}.$$ So,
$$J=\|\nabla U^{(j_0)}\|_{L^2}^2=\frac{\|Q\|_{L^2}^2}{2}.$$

So far, we have obtained the existence of the positive solution of
(\ref{ground}). In addition, Theorem 3 in \cite{KrLR08} together
with Theorem 1.2 in \cite{Liu} implies that this positive solution
is also radial and unique in $H^1(\mathbb{R}^4)$. Note that the
uniqueness proof strongly depends on the specific features of
equation $(\ref{ground})$. In fact, the uniqueness of the ground
state $Q$ of $(\ref{ground})$ has not be resolved completely for the
general potential $V(x)$, and be stated as an open problem in
\cite{FrL04}.

We first make use of the ground state $Q$ to give a sufficient
condition for the global existence of (\ref{equ42}), which together
with $(\ref{pcs})$ implies that $\big\|Q\big\|_{L^2}$ is the minimal
mass of the blow up solutions.

\begin{theorem}\label{thm42}
If $u_0\in H^1(\mathbb{R}^4)$ and $\|u_0\|_{L^2}<\|Q\|_{L^2}$, then
the solution $u(t)$ of (\ref{equ42}) is global in time.
\end{theorem}
{\noindent  \it Proof.} By the local wellposedness theory, it
suffices to prove that for every $t\in \Bbb R$, we have
$$\|\nabla u(t)\|_{L^2}<+\infty.$$
Now from Proposition $\ref{pro41}$ and the conservation of mass,  we
have
\begin{align}\label{equ61}
E(u(t))&=\frac{1}{2}\int|\nabla
u(t)|^2dx-\frac{1}{4}\int(V*|u(t)|^2)|u(t)|^2dx\nonumber\\
&\geq\frac{1}{2}\|\nabla
u(t)\|_{L^2}^2-\frac{1}{4}\frac{2}{\|Q\|_{L^2}^2}\|u(t)\|_{L^2}^2\|\nabla
u(t)\|_{L^2}^2\nonumber\\
&=\frac{1}{2}\|\nabla
u(t)\|_{L^2}^2\Big(1-\frac{\|u_0\|_{L^2}^2}{\|Q\|_{L^2}^2}\Big).
\end{align}
Since $\|u_0\|_{L^2}<\|Q\|_{L^2}$, so we have the uniform bound of
$\|\nabla u(t)\|_{L^2}^2$. This proves the global existence.

Before we prove Theorem \ref{thm4}, we state a proposition in two
equivalent forms.
\begin{proposition}[Static version]\label{pro42}
If $u\in H^1(\mathbb{R}^4)$ such that $\|u\|_{L^2}=\|Q\|_{L^2}$ and
$E(u)=0$, then $u(x)$ is of the following form
$$u(x)=e^{i\theta}\lambda^2Q(\lambda x+b),
\ \ \text{for some}\ \ \theta\in\mathbb{R},\ \lambda>0, \
b\in\mathbb{R}^4.$$
\end{proposition}

{\noindent \it Proof.} Since $E(u)=0$, we have $\|\nabla
u\|_{L^2}^2=\frac{1}{2}\|u\|_{L^V}^4$. So we get
$$J(u)=\frac{\|Q\|_{L^2}^2\|\nabla u\|_{L^2}^2}{\|u\|_{L^V}^4}=\frac{1}{2}\|Q\|_{L^2}^2=J.$$
By Proposition \ref{pro41} and the uniqueness of the ground state
$Q$, $u$ is of the form $u(x)=aQ(\lambda x+b)$. The condition
$\|u\|_{L^2}=\|Q\|_{L^2}$ ensures that $|a|=\lambda^2$. So
$u(x)=e^{i\theta}\lambda^2Q(\lambda x+b)$.

\begin{proposition}[Dynamic version]\label{pro43}
Let $\{u_n\}_{n=1}^\infty$ be a sequence in $H^1(\mathbb{R}^4)$ such
that $\|u_n\|_{L^2}=\|Q\|_{L^2}$, $E(u_n)\leq M$ and $\|\nabla
u_n\|_{L^2}\rightarrow\infty$. We define
$$\lambda_n:=\frac{\|\nabla u_n\|_{L^2}}{\|\nabla Q\|_{L^2}},$$
then there exists a subsequence {\rm(}still denoted by
$\{u_n\}${\rm)}, a sequence $(y_n)\subset\mathbb{R}^4$ and a real
number $\theta$ such that
\begin{equation}\label{equ66}
e^{i\theta}\lambda_n^{-2}u_n(\lambda_n^{-1}x+y_n)\rightarrow
Q(x)\ \text{strongly in}\ H^1.
\end{equation}
\end{proposition}
{\noindent \it Proof.} Let
$$\tilde{u}_n(x)=\frac{1}{\lambda_n^{2}}u_n(\frac{x}{\lambda_n}),$$
then $\|\tilde{u}_n\|_{L^2}=\|Q\|_{L^2}$ and
$\|\nabla \tilde{u}_n\|_{L^2}=\|\nabla Q\|_{L^2}$. Moreover,
$$E(\tilde{u}_n)=\frac{E(u_n)}{\lambda_n^2}\rightarrow0,\ \ \text{as}\ \ n\rightarrow\infty.$$
So we have
\begin{align*}
J(\tilde{u}_n)&=\|Q\|_{L^2}^{2}\frac{\|\nabla\tilde{u}_n\|_{L^2}^2}{\|\tilde{u}_n\|_{L^V}^4}
=\|Q\|_{L^2}^{2}\frac{\|\nabla\tilde{u}_n\|_{L^2}^2}{2\|\nabla\tilde{u}_n\|_{L^2}^2-4E(\tilde{u}_n)}
\ \longrightarrow\ \frac{\|Q\|_{L^2}^{2}}{2}=J, \quad
n\longrightarrow\infty.
\end{align*}
Therefore, by Lemma \ref{lem43}, we can choose a subsequence
$\tilde{u}_n $ and $(x_n)\subset\mathbb{R}^4$ such that
$\tilde{u}_n(x+x_n)\rightarrow aQ(\lambda x+b)$ in $H^1$. The
conditions $\|\tilde{u}_n\|_{L^2}=\|Q\|_{L^2}$ and $\|\nabla
\tilde{u}_n\|_{L^2}=\|\nabla Q\|_{L^2}$ imply $|a|=\lambda=1$, so we
have (\ref{equ66}) for
$y_n=\lambda_n^{-1}(x_n-b)$.\\

In order to prove Theorem \ref{thm4}, we also need the following
lemma. The proof relies heavily on the techniques in V. Banica
\cite{Ba}.
\begin{lemma}\label{lem45}
Suppose $u\in H^1(\mathbb{R}^4)$,
$\|u\|_{L^2}=\|Q\|_{L^2}$, then for all real function $w\in C^1$
with $\nabla w$ is bounded, we have$$\Big|\int_{\mathbb{R}^4}\nabla
w(x)\Im(u\nabla u)(x)dx\Big|\leq\sqrt{2}
E(u)^{\frac{1}{2}}\Big(\int|u|^2|\nabla
w|^2dx\Big)^{\frac{1}{2}}.$$
\end{lemma}

{\noindent  \it Proof.} Since
$$\|ue^{isw(x)}\|_{L^2}=\|u\|_{L^2}=\|Q\|_{L^2},$$ for any
$s\in\mathbb{R}$, by (\ref{equ61}) we know that
$E(ue^{isw(x)})\geq0$. So, for any $s$,
$$\frac{1}{2}\int_{\mathbb{R}^4}|\nabla u+isu\nabla w|^2dx-\frac{1}{4}\int_{\mathbb{R}^4}(V*|u|^2)|u|^2dx\geq0.$$
Namely, $$E(u)+s\int_{\mathbb{R}^4}\nabla w\Im(u\nabla
u)dx+\frac{s^2}{2}\int_{\mathbb{R}^4}|u|^2|\nabla w|^2dx\geq0.$$
Note that this holds for any $s$, so the discriminant is
non-positive. So we get the result.

Now we turn to the proof of Theorem \ref{thm4} and Theorem
\ref{thm5}, which is borrowed from \cite{HMKer}.

\vskip0.2cm
 {\noindent  \it Proof of Theorem \ref{thm4}.} Suppose
$u(t,x)$ is the solution of (\ref{equ42}) which blows up at $T$ and
let $\{t_n\}_{n=1}^\infty$ be an arbitrary sequence such that
$t_n\uparrow T$. Let $u_n=u(t_n)$, by Proposition \ref{pro43}, we
have
$$e^{i\theta}\lambda_n^{-2}u_n(\lambda_n^{-1}x+y_n)\rightarrow
Q(x)\ \ \text{strongly in}\ H^1.$$ From this we get
\begin{equation}\label{equ62}|u(t_n,x)|^2dx-\|Q\|_{L^2}^2\delta_{x=y_n}\rightharpoonup0.\end{equation}where
$y_n\rightarrow0$ (up to translation) or $y_n\rightarrow\infty$.

Now let $\phi\in C_0^\infty(\mathbb{R}^4)$ be a nonnegative radial
function such that $$\phi(x)=|x|^2,\ \ \text{if}\ \ |x|<1\
\text{and}\ \ |\nabla \phi|^2\leq C\phi(x).$$ For every
$p\in\mathbb{N}^*$ we define
$$\phi_p(x)=p^2\phi(\frac{x}{p})\ \ \text{and}\ \ g_p(t)=\int\phi_p(x)|u(t,x)|^2dx.$$
By Lemma \ref{lem45}, for every $t\in [0,T)$, we have\begin{align*}
|\dot{g}_p(t)|&=2\Big|\int_{\mathbb{R}^4}\nabla \phi_p(x)\Im(u\nabla
u)(x)dx\Big|\leq 2\sqrt{2}E(u_0)^{\frac{1}{2}}\Big(\int|u|^2|\nabla
\phi_p(x)|^2dx\Big)^{\frac{1}{2}}\\
&\leq
CE(u_0)^{\frac{1}{2}}\Big(\int|u|^2\phi_p(x)dx\Big)^{\frac{1}{2}}\leq
C(u_0)\sqrt{g_p(t)}.\end{align*} Integrating with respect to $t$, we
get that $$\Big|\sqrt{g_p(t)}-\sqrt{g_p(t_n)}\Big|\leq
C(u_0)|t_n-t|.$$ If $y_n\rightarrow0$, then
$g_p(t_n)\rightarrow\|Q\|_{L^2}^2\phi_p(0)=0$ by (\ref{equ62}); if
$|y_n|\rightarrow\infty$, also $g_p(t_n)\rightarrow0$ since $\phi_p$
is compactly supported. So, if we let $n$ go to infinity, we have
$$g_p(t)\leq C(u_0)(T-t)^2.$$

Now fix $t\in[0,T)$ and let $p$ go to infinity, then by
(\ref{equ21}) we get
\begin{equation}\label{equ63}8t^2E(e^{i\frac{|x|^2}{4t}}u_0)=\int|x|^2|u(t,x)|^2dx\leq
C(u_0)(T-t)^2.\end{equation} Hence $$|y_n|^2\|Q\|_{L^2}^2\leq
C(u_0)T^2.$$ Thus $y_n$ can not go to infinity. This implies that
$\{y_n\}$ converges to 0. Let $t$ goes to $T$, from (\ref{equ63}),
we get
$$E(e^{i\frac{|x|^2}{4T}}u_0)=0.$$ Note also that
$$\|e^{i\frac{|x|^2}{4T}}u_0\|_{L^2}=\|Q\|_{L^2}.$$ By
Proposition \ref{pro42}, we conclude that
$e^{i\frac{|x|^2}{4T}}u_0\in\mathcal{A}$.

\vskip0.2cm
 {\noindent \it Proof of Theorem \ref{thm5}.}  We denote
$$\rho(t)=\frac{\|\nabla Q\|_{L^2}}{\|\nabla u\|_{L^2}}\ \
\text{and}\ \ v(t,x)=\rho^2u(t,\rho x).$$ Let $\{t_n\}_{n=1}^\infty$
be an arbitrary time sequence such that $t_n\uparrow T$,
$v_n(x)=v(t_n,x)$, then by mass conservation and the definition of
$\rho(t)$, we have $$\|v_n\|_{L^2}=\|u_0\|_{L^2}\ \ \text{and}\ \
\|\nabla v_n\|_{L^2}=\|\nabla Q\|_{L^2}.$$ Since $u$ blows up at
time $T$, we have
$$\rho(t_n)\rightarrow0, \ \ \text{as} \ t_n\rightarrow T.$$
So we have $$E(v_n)=\rho_n^2E(u_0)\rightarrow0, \text{as}\
n\rightarrow\infty.$$ In particular,
$$\|v_n\|_{L^V}^4\rightarrow2\|\nabla Q\|_{L^2}^2, \text{as}\ n\rightarrow\infty.$$

According to Lemma \ref{lem43}, the sequence $\{v_n\}_{n=1}^\infty$
can be written, up to a subsequence, as
$$v_n(x)=\sum_{j=1}^lU^{(j)}(x-x_n^j)+r_n^l(x)$$
such that (\ref{equ45}), (\ref{equ46}) and (\ref{equ47}) hold. This
implies, in particular, that
$$2\|\nabla Q\|_{L^2}^2\leq\limsup_{n\rightarrow\infty}\|v_n\|_{L^V}^4
=\limsup_{n\rightarrow\infty}\Big\|\sum_{j=1}^\infty
U^j(\cdot-x_n^j)\Big\|_{L^V}^4.$$ As in the discussion of the proof
of Proposition \ref{pro41}, the pairwise orthogonality of the family
$\{x^j\}_{j=1}^\infty$, together with (\ref{ground}) and
(\ref{equ47}), gives
\begin{align*}
2\|\nabla Q\|_{L^2}^2  \leq & \sum_{j=1}^\infty \|U^j\|_{L^V}^4\leq
\sum_{j=1}^\infty\frac{2}{\|Q\|_{L^2}^2}\|U^j\|_{L^2}^2\|\nabla
U^j\|_{L^2}^2\\
\leq&\frac{2}{\|Q\|_{L^2}^2}\sup_{j\geq1}\|U^j\|_{L^2}^2\sum_{j=1}^\infty\|\nabla
U^j\|_{L^2}^2  \leq \frac{2}{\|Q\|_{L^2}^2}\|\nabla
v_n\|_{L^2}^2\sup_{j\geq1}\|U^j\|_{L^2}^2\\
\leq&\frac{2}{\|Q\|_{L^2}^2}\|\nabla
Q\|_{L^2}^2\sup_{j\geq1}\|U^j\|_{L^2}^2.
\end{align*}
Therefore, we get that
$$\sup_{j\geq1}\|U^j\|_{L^2}^2\geq\|Q\|_{L^2}^2.$$
Since $\sum \|U^j\|_{L^2}^2$ converges, the supremum above is
attained. In particular, there exists $j_0$ such that
$$\|U^{j_0}\|_{L^2}^2\geq\|Q\|_{L^2}^2.$$
On the other hand, a change of variables gives
$$v_n(x+x_n^{j_0})=U^{j_0}(x)+\sum_{1\leq j\leq l \atop j\ne j_0}U^j(x+x_n^{j_0}-x_n^j)+\tilde{r}_n^l(x),$$
where $\tilde{r}_n^l(x)=r_n^l(x+x_n^{j_0})$. The pairwise
orthogonality of the family $\{x^j\}_{j=1}^\infty$ implies
$$U^j(\cdot+x_n^{j_0}-x_n^j)\rightharpoonup0, \ \text{weakly}$$
for every $j\ne j_0$. Hence we get
$$r_n(\cdot+x_n^{j_0})\rightharpoonup U^{j_0}+\tilde{r}^l,$$
where $\tilde{r}^l$ denote the weak limit of
$\{\tilde{r}_n^l\}_{n=1}^\infty$. However, we have
\begin{equation*}
\|\tilde{r}^l\|_{L^V}\leq
\limsup_{n\rightarrow\infty}\|\tilde{r}_n^l\|_{L^V}=\limsup_{n\rightarrow\infty}
\|r_n^l\|_{L^V}\buildrel{l\rightarrow\infty}\over{\longrightarrow}0.
\end{equation*}
By uniqueness of weak limit, we get $$\tilde{r}^l=0$$ for every
$l\ne j_0$ so that $$r_n(\cdot+x_n^{j_0})\rightharpoonup U^{j_0}, \
\text{in}\ H^1,$$ namely, $$\rho_n^2u(t_n,
\rho_n\cdot+x_n^{j_0})\rightharpoonup U^{j_0}\in H^1\ \
\text{weakly}.$$ Thus for every $A>0$,
\begin{equation*}
\liminf_{n\rightarrow+\infty}\int_{|x|\leq A}\rho_n^4|u(t_n,
\rho_nx+x_n)|^2dx\geq\int_{|x|\leq A}|U^{j_0}|^2dx.
\end{equation*}
In view of the assumption $\lambda(t_n)/\rho_n\rightarrow\infty$,
this gives immediately
\begin{equation*}
\liminf_{n\rightarrow+\infty}\sup_{y\in\Bbb R^4}\int_{|x-y|\leq
\lambda(t_n)}|u(t_n,x)|^2dx\geq\int_{|x|\leq A}|U^{j_0}|^2dx
\end{equation*}
for every $A>0$, which means that
\begin{equation*}
\liminf_{n\rightarrow+\infty}\sup_{y\in\Bbb R^4}\int_{|x-y|\leq
\lambda(t_n)}|u(t_n,x)|^2dx\geq\int|U^{j_0}|^2dx\geq\int|Q|^2dx.
\end{equation*}
Since the sequence $\{t_n\}_{n=1}^\infty$ is arbitrary, we infer
\begin{equation*}
\liminf_{t\rightarrow T}\sup_{y\in \Bbb
R^4}\int_{|x-y|\leq\lambda(t)}|u(t,x)|^2dx\geq\int|Q|^2dx.
\end{equation*}
But for every $t\in [0,T)$, the function $y\mapsto
\int_{|x-y|\leq\lambda(t)}|u(t,x)|^2dx$ is continuous and goes to 0
at infinity. As a result, we get
\begin{equation*}
\sup_{y\in\Bbb
R^4}\int_{|x-y|\leq\lambda(t)}|u(t,x)|^2dx=\int_{|x-x(t)|\leq\lambda(t)}|u(t,x)|^2dx,
\end{equation*}
for some $x(t)\in \Bbb R^4$ and Theorem \ref{thm5} is proved.

  \section{The blow-up dynamics of the focusing mass critical Hartree equation with $L^2$ data}
\setcounter{equation}{0} In this section we prove Theorem
\ref{thm-mass1} and Theorem \ref{thm-mass2}.
\begin{definition}\label{def1}
For every sequence $\mathbf{\Gamma}_n=\{\rho_n, t_n, \xi_n,
x_n\}_{n=1}^\infty\subset
\mathbb{R}_+^*\times\mathbb{R}\times\mathbb{R}^4\times\mathbb{R}^4$,
we define the isometric operator $\mathbf{\Gamma}_n$ on
$L^3_{t,x}(\mathbb{R}\times\mathbb{R}^4)$ by
$$\mathbf{\Gamma}_n(f)(t,x)
=\rho_n^2e^{ix\cdot\xi_n}e^{-it|\xi_n|^2}f(\rho_n^2t+t_n,
\rho_n(x-t\xi_n)+x_n).$$ Two sequences
$\mathbf{\Gamma}^j=\{\rho_n^j,t_n^j,\xi_n^j,x_n^j\}_{n=1}^{\infty}$
and
$\mathbf{\Gamma}^k=\{\rho_n^k,t_n^k,\xi_n^k,x_n^k\}_{n=1}^{\infty}$
are said to be orthogonal if
$$\frac{\rho_n^j}{\rho_n^k}+\frac{\rho_n^k}{\rho_n^j}\rightarrow+\infty $$
or
$$\rho_n^j=\rho_n^k\ \ and\ \ \frac{|\xi_n^j-\xi_n^k|}{\rho_n^j}+|t_n^j-t_n^k|
+\bigg|\frac{\xi_n^j-\xi_n^k}{\rho_n^j}t_n^j+x_n^j-x_n^k\bigg|\rightarrow+\infty.$$
\end{definition}
\begin{lemma}[Linear profile decomposition \cite{BeV}] \label{lem44}Let $\{\varphi_n\}_{n\in\mathbb{N}}$
be a bounded sequence in $L^2(\mathbb{R}^4)$. Then there exists a
subsequence of $\{\varphi_n\}_{n=1}^\infty$ (still denoted by
$\{\varphi_n\}_{n=1}^\infty$) which satisfies the following
properties: there exists a family $\{V^j\}_{j=1}^{\infty}$ of
solutions of (\ref{equ41}) and a family of pairwise orthogonal
sequences
$\mathbf{\Gamma}^j=\{\rho_n^j,t_n^j,\xi_n^j,x_n^j\}_{n=1}^{\infty}$,
 such that for every $(t,x)\in\mathbb{R}\times\mathbb{R}^4$, we have
\begin{equation}\label{p1}e^{it\Delta}\varphi_n(x)=\sum_{j=1}^l\mathbf{\Gamma}_n^jV^j(t,x)+w_n^l(t,x),\end{equation}
with
\begin{equation}\label{p2}\limsup_{n\rightarrow\infty}\|w_n^l\|_{L^3
(\mathbb{R}\times\mathbb{R}^4)}\rightarrow0,\ \text{as}\
l\rightarrow\infty.\end{equation} Moreover, for every $l\geq1$,
\begin{equation}\label{p3}\|\varphi_n\|_{L^2}^2=\sum_{j=1}^l\|V^j\|_{L^2}^2+\|w_n^l\|_{L^2}^2+o_n(1).\end{equation}\end{lemma}

\begin{definition}\label{def2}
Let $\Gamma_n=\{\rho_n, t_n, \xi_n, x_n\}_{n=1}^\infty$ be a
sequence of $\mathbb{R}^*_+\times
\mathbb{R}\times\mathbb{R}^4\times\mathbb{R}^4$  such that the
quantity $\{t_n\}_{n=1}^\infty$ has a limit in $[-\infty,+\infty]$
when $n$ goes to the infinity. Let $V$ be a solution of linear
Schr\"{o}dinger equation (\ref{equ41}). We say that $U$ is the
nonlinear profile associated to $\{V,\Gamma_n\}_{n=1}^\infty$ if $U$
is the unique maximal solution of the nonlinear Schr\"{o}dinger
equation (\ref{equ42}) satisfying
$$\Big\|(U-V)(t_n, \cdot)\Big\|_{L^2(\mathbb{R}^4)}\rightarrow0,\ \ \text{as} \ n\rightarrow\infty.$$
\end{definition}

In order to prove Theorem \ref{thm-mass1} and Theorem
\ref{thm-mass2}, we first state a key theorem, which is similar to
that in \cite{Ker1} and \cite{Ker2} and its proof is the same
essence with that of stability theory.

\begin{theorem}[Nonlinear profile decomposition]\label{thm41} Let $\{\varphi_n\}_{n=1}^\infty$ be a
bounded family of $L^2(\mathbb{R}^4)$ and $\{u_n\}_{n=1}^\infty$ the
corresponding family of solutions to (\ref{equ42}) with initial data
$\{\varphi_n\}_{n=1}^\infty$. Let $\{V^j,
\mathbf{\Gamma}_n^j\}_{j=1}^\infty$ be the family of linear profiles
associated to $\{\varphi_n\}_{j=1}^\infty$ via Lemma \ref{lem44} and
$\{U^j\}_{j=1}^\infty$ the family of nonlinear profiles associated
to $\{V^j, \mathbf{\Gamma}_n^j\}_{j=1}^\infty$ via Definition
\ref{def2}. Let $\{I_n\}_{n=1}^\infty$ be a family of intervals
containing the
origin 0. Then the following statements are equivalent:\\
(i) For every $j\geq1$, we have
$$\lim_{n\rightarrow\infty}\|\mathbf{\Gamma}_n^jU^j\|_{L_{t,x}^3[I_n]}<\infty,$$
(ii)
$$\lim_{n\rightarrow\infty}\|u_n\|_{L_{t,x}^3[I_n]}<\infty.$$
Moreover, if (i) or (ii) holds, then
\begin{equation}\label{equ22}u_n=\sum_{j=1}^l\mathbf{\Gamma}_n^jU^j+w_n^l+r_n^l,\end{equation} where $w_n^l$ is
as in (\ref{p2}) and
\begin{equation}\label{equ23}
\lim_{n\rightarrow\infty}(\|r_n^l\|_{L_{t,x}^3[I_n]}+\sup_{t\in
I_n}\|r_n^l\|_{L^2})\rightarrow0 \quad \text{as}\quad
l\rightarrow\infty
\end{equation}
\end{theorem}

\noindent{\it Proof.}  {\bf Step 1:} We prove (\ref{equ22}) and
(\ref{equ23}) provided that (i) or (ii) holds.
 Let
$$r_n^l=u_n-\sum_{j=1}^lU_n^j-w_n^l,\  \text{where} \ \ U_n^j:=\mathbf{\Gamma}_n^jU^{j},$$
 and let $V_n^j:=\mathbf{\Gamma}_n^jV^{j}$, then $r_n^l$ satisfies the
following equation
\begin{equation}
\left\{ \aligned
    i\partial_tr_n^l +  \Delta r_n^l  & = f_n^l, \\
     r_n^l(0)&=\sum_{j=1}^l(V_n^j-U_n^j)(0,x).
\endaligned
\right.
\end{equation}
where
$$f_n^l:=p(W_n^l+w_n^l+r_n^l)-\sum_{j=1}^lp(U_n^j),$$
and $$p(z):=-(|x|^{-2}*|z|^2)z,\ \ W_n^l:=\sum_{j=1}^lU_n^j.$$
 It suffices to prove that
\begin{equation}\label{equ64}
\lim_{n\rightarrow\infty}(\|r_n^l\|_{L_{t,x}^3[I_n]}+\sup_{t\in
I_n}\|r_n^l\|_{L^2})\buildrel{l\rightarrow\infty}\over\longrightarrow0.
\end{equation}

By Strichartz estimates and Young's inequality, we have
\begin{align}
\big\|r_n^l\big\|_{L_{t,x}^3[I_n]}+\sup_{t\in I_n}\|r_n^l\|_{L^2}
\lesssim&\quad
\Big\|p(W_n^l+w_n^l+r_n^l)-\sum_{j=1}^lp(U_n^j)\Big\|_{\dot{\mathcal {N}}^0[I_n]}+\|r_n^l(0,\cdot)\|_{L^2}\nonumber\\
\lesssim& \quad \Big\|p(W_n^l)-\sum_{j=1}^lp(U_n^j)\Big\|_{\dot{\mathcal {N}}^0[I_n]}\label{equ411}\\
&+\Big\|p(W_n^l+w_n^l)-p(W_n^l)\Big\|_{L^1_tL^2_x[I_n]}\label{equ412}\\
&+\Big\|p(W_n^l+w_n^l+r_n^l)-p(W_n^l+w_n^l)\Big\|_{L^1_tL^2_x[I_n]}\label{equ413}\\
&+\|r_n^l(0,\cdot)\|_{L^2}.\nonumber
\end{align}
We will estimate these three terms, respectively. Firstly, we
estimate (\ref{equ411}).
\begin{align}
(\ref{equ411})\leq& \sum_{j_1=1}^l\sum_{j_2\neq
j_1}\Big\|(|x|^{-2}*|U_n^{j_1}|^2)U_n^{j_2}\Big\|_{L_{t,x}^\frac{3}{2}[I_n]}\label{equ414}\\
&+\sum_{j_1=1}^l\sum_{j_2\neq
j_1}\sum_{j_3=1}^l\Big\|\big(|x|^{-2}*(U_n^{j_1}U_n^{j_2})\big)U_n^{j_3}\Big\|_{L^1_tL^2_x[I_n]}.\label{equ415}
\end{align}
Without loss of generality we can assume that both $U^{j_1}$ and
$U^{j_2}$ have compact support in $t$ and $x$. Let $V(x)=|x|^{-2}$,
then we have
\begin{align*}
&\iint |(V*|U_n^{j_1}|^2)U_n^{j_2}|^{\frac{3}{2}}dxdt\\
=&\iint
\Big|\int(\rho_n^{j_1})^4|U^{j_1}((\rho_n^{j_1})^2t+t_n^{j_1},
\rho_n^{j_1}(x-y-t\xi_n^{j_1})+x_n^{j_1})|^2V(y)dy\\
&\qquad \qquad \qquad \qquad \qquad \times
(\rho_n^{j_2})^2U^{j_2}((\rho_n^{j_2})^2t+t_n^{j_2},
\rho_n^{j_2}(x-t\xi_n^{j_2})+x_n^{j_2})\Big|^{\frac{3}{2}}dxdt\\
=&\bigg(\frac{\rho_n^{j_2}}{\rho_n^{j_1}}\bigg)^3\iint
\Bigg|\int|U^{j_1}
(\tilde{t},\tilde{x}-\tilde{y})|^2V(\tilde{y})d\tilde{y}U^{j_2}\Bigg(\bigg(\frac{\rho_n^{j_2}}{\rho_n^{j_1}}\bigg)^2\tilde{t}-
\bigg(\frac{\rho_n^{j_2}}{\rho_n^{j_1}}\bigg)^2t_n^{j_1}+t_n^{j_2},\\
&\qquad \qquad \qquad
\frac{\rho_n^{j_2}}{\rho_n^{j_1}}\tilde{x}+\frac{\rho_n^{j_2}(\xi_n^1-\xi_n^2)}{(\rho_n^{j_1})^2}\tilde{t}
-\frac{\rho_n^{j_2}(\xi_n^{j_1}-\xi_n^{j_2})}{(\rho_n^{j_1})^2}t_n^{j_1}
-\frac{\rho_n^{j_2}x_n^{j_1}}{\rho_n^{j_1}}+x_n^{j_2}\Bigg)\Bigg|^{\frac{3}{2}}d\tilde{x}d\tilde{t}.
\end{align*}

If
$\rho_n^{j_2}/\rho_n^{j_1}+\rho_n^{j_1}/\rho_n^{j_2}\rightarrow+\infty$
or $|t_n^{j_1}-t_n^{j_2}|\rightarrow+\infty$,
 by the compact support assumption on $t$, we conclude that
 $(\ref{equ414})\rightarrow0$.
 Otherwise, by orthogonality we have
\begin{equation}\label{equ422}\frac{|\xi_n^{j_1}-\xi_n^{j_2}|}{\rho_n^{j_1}}+\bigg|\frac{\xi_n^{j_1}-\xi_n^{j_2}}{\rho_n^{j_1}}t_n^{j_1}+x_n^{j_1}-x_n^{j_2}\bigg|\rightarrow+\infty.\end{equation}
Without loss of generality, we may assume that
$\rho_n^{j_2}/\rho_n^{j_1}\rightarrow1$. Then the complicated
expression of the function $U^{j_2}$ of $\tilde{t}$ and $\tilde{x}$
can be simplified to
$$U^{j_2}\bigg(\tilde{t}-t_n^{j_1}+t_n^{j_2},
\frac{\xi_n^{j_1}-\xi_n^{j_2}}{\rho_n^{j_1}}\tilde{t}+\tilde{x}-x_n^{j_1}+x_n^{j_2}-\frac{\xi_n^{j_1}-\xi_n^{j_2}}{\rho_n^{j_1}}t_n^{j_1}\bigg).$$
Meanwhile, we have $$\int|U^{j_1}
(\tilde{t},\tilde{x}-\tilde{y})|^2V(\tilde{y})d\tilde{y}\leq
\int_{|\tilde{y}|\leq1}|U^{j_1}
(\tilde{t},\tilde{x}-\tilde{y})|^2V(\tilde{y})d\tilde{y}+\sum_{j=0}^\infty\int_{2^j\leq|\tilde{y}|\leq
2^{j+1}}|U^{j_1}
(\tilde{t},\tilde{x}-\tilde{y})|^2V(\tilde{y})d\tilde{y}.$$ Note
that $U^{j_1}$ is compactly supported in $x$, so for any fixed $j$,
$$\int_{2^j\leq|\tilde{y}|\leq 2^{j+1}}|U^{j_1}
(\tilde{t},\cdot-\tilde{y})|^2V(\tilde{y})d\tilde{y}$$ is also
compactly supported. Thus (\ref{equ422}) implies that for any
$j_1\neq j_2$,
\begin{align*}
\lim_{n\rightarrow\infty}\iint \Bigg|\int_{2^j\leq|\tilde{y}|\leq
2^{j+1}}& |U^{j_1}
(\tilde{t},\cdot-\tilde{y})|^2V(\tilde{y})d\tilde{y}U^{j_2}\Bigg(\tilde{t}-t_n^{j_1}+t_n^{j_2},\\
&
\frac{\xi_n^{j_1}-\xi_n^{j_2}}{\rho_n^{j_1}}\tilde{t}+\tilde{x}-x_n^{j_1}+x_n^{j_2}
-\frac{\xi_n^{j_1}-\xi_n^{j_2}}{\rho_n^{j_1}}t_n^{j_1}\Bigg)\Bigg|^{\frac{3}{2}}d\tilde{x}d\tilde{t}
=0.\end{align*} Therefore, we get that $(\ref{equ414})\rightarrow0$
as $n\rightarrow\infty$.

On the other hand,
$$\Big\|(|x|^{-2}*(U_n^{j_1}U_n^{j_2})U_n^{j_3}\Big\|_{L^1_tL^2_x[I_n]}\leq
C\Big\|U^{j_1}_nU^{j_2}_n\Big\|_{L_{t,x}^\frac{3}{2}}\Big\|U^{j_3}_n\Big\|_{L^3_{t,x}}.$$
By orthogonality,
$$\Big\|U^{j_1}_nU^{j_2}_n\Big\|_{L_{t,x}^\frac{3}{2}}\rightarrow0,\ \ \text{as}\ \ n\rightarrow\infty.$$
Because $\Big\|U^{j_3}_n\Big\|_{L_{t,x}^3}$ is bounded, we have
$$(\ref{equ415})\xrightarrow{n\rightarrow\infty}0.$$

Next, we prove that
$$\lim_{l\rightarrow\infty}\Big(\lim_{n\rightarrow\infty}\big\|W_n^l+w_n^l\big\|_{L_{t,x}^3[I_n]}\Big)\leq C.$$
From (\ref{p3}), we have
$$\big\|w_n^l\big\|_{L_{t,x}^3[I_n]}\leq
C\|w_n^l(0)\|_{L^2}\leq C\|\varphi_n\|_{L^2}.$$ It suffices to
verify
\begin{equation}\label{equ420}\lim_{l\rightarrow\infty}\Big(\lim_{n\rightarrow\infty}\big\|W_n^l\big\|_{L_{t,x}^3[I_n]}\Big)\leq C.\end{equation} From the orthogonality of $\Gamma_n^j$, as in
\cite{Ker1}, we can get that for every $l\geq1$
\begin{align*}
\big\|W_n^l\big\|_{L_{t,x}^3[I_n]}^3
=\big\|\sum_{j=1}^lU_n^j\big\|_{L_{t,x}^3[I_n]}^3
\rightarrow\sum_{j=1}^l\big\|U^j_n\big\|_{L_{t,x}^3[I_n]}^3,\ \
\text{as}\ \ n\rightarrow\infty.
\end{align*}
Meanwhile by (\ref{p3}), the series $\sum\|V^j\|_{L^2}^2$ converge.
Thus for every $\epsilon>0$, there exists $l(\epsilon)$ such that
$$\|V^j\|_{L^2}\leq\epsilon, \quad, \forall j>l(\epsilon).$$
The theory of small data asserts that , for $\epsilon$ sufficiently
small, $U^{j}$ is global and
$$\|U^j\|_{L_{t,x}^3}\lesssim\|V^j\|_{L^2},$$
which yields that
$$\sum_{j>l(\epsilon)}\|U^j\|_{L_{t,x}^3}^3<\infty.$$
So we have to deal only with a finite number of nonlinear profiles
$\{U^{j}\}_{1\leq j\leq l(\epsilon)}$. But in view of the pairwise
orthogonality of $\{\mathbf{\Gamma}_n^j\}_{j=1}^\infty$, one has
$$\lim_{n\rightarrow\infty}\Big\|\sum_{j=1}^{l(\epsilon)}U_n^j\Big\|_{L_{t,x}^3[I_n]}
\leq\sum_{j=1}^{l(\epsilon)}\lim_{n\rightarrow\infty}\Big\|U_n^j\Big\|_{L_{t,x}^3[I_n]}<\infty$$
and then (\ref{equ420}) follows.

Now, we estimate (\ref{equ412}).
\begin{align*}
&\Big\|p(W_n^l+w_n^l)-p(W_n^l)\Big\|_{L^1_tL^2_x[I_n]}\\
\lesssim &
 \Big\|\big(|x|^{-2}*|W_n^l+w_n^l|^2\big)w_n^l\Big\|_{L^1_tL^2_x[I_n]}
+\Big\|\big(|x|^{-2}*(W_n^lw_n^l)\big)w_n^l\Big\|_{L^1_tL^2_x[I_n]} +\Big\|\big(|x|^{-2}*|w_n^l|^2\big)W_n^l\Big\|_{L^1_tL^2_x[I_n]} \\
\lesssim&\Big\|W_n^l\Big\|_{L_{t,x}^3[I_n]}^2
\big\|w_n^l\big\|_{L_{t,x}^3[I_n]}+
\big\|w_n^l\big\|^2_{L_{t,x}^3[I_n]}
\Big(\big\|W_n^l\big\|_{L_{t,x}^3[I_n]}+\big\|w_n^l\big\|_{L_{t,x}^3[I_n]}\Big)\\
=&o_n(1).
\end{align*}
The last equality is due to (\ref{equ420}) and the fact that
$\big\|w_n^l\big\|_{L_{t,x}^3[I_n]}\rightarrow0$ as
$l\rightarrow\infty$.

(\ref{equ413}) can be estimated similarly. In fact, we have
\begin{align*}
(\ref{equ413})\lesssim\Big(\|W_n^l+w_n^l\|_{L_{t,x}^3[I_n]}^2\|r_n^l\|_{L_{t,x}^3[I_n]}
+\|W_n^l+w_n^l\|_{L_{t,x}^3[I_n]}\|r_n^l\|_{L_{t,x}^3[I_n]}^2+\|r_n^l\|_{L_{t,x}^3[I_n]}^3\Big).
\end{align*}

Now we can prove (\ref{equ64}). Collecting all the previous facts,
we have
\begin{align}\label{equ416}&\sup_{t\in I_n}\|r_n^l\|_{L^2}+\|r_n^l\|_{L_{t,x}^3[I_n]}\nonumber\\
\leq
&C\Big(\|W_n^l+w_n^l\|_{L_{t,x}^3[I_n]}\|r_n^l\|_{L_{t,x}^3[I_n]}
+\|r_n^l\|_{L_{t,x}^3[I_n]}^3+\|r_n^l\|_{L_{t,x}^3[I_n]}^2+\|r_n^l(0,\cdot)\|_{L^2}\Big)+o_n(1).\end{align}
As in \cite{Ker2}, for every $\varepsilon>0$ we can divide
$I_n^+=I_n\cap\mathbb{R}_+$ into finite n-dependent intervals,
namely,
$$I_n^+=[0, a_n^1]\cup[a_n^1,a_n^2]\cup\cdots\cup[a_n^{p-1},a_n^p),$$ with each interval
 denoted by $I_n^i\  (i=1,2,\cdots,p)$, such that for every $1\leq i\leq p$ and every $l\geq1$,
$$\limsup_{n\rightarrow\infty}\|W_n^l+w_n^l\|_{L_{t,x}^3(I_n^i\times\mathbb{R}^4)}\leq\varepsilon.$$
The $I_n^-=I_n\cap\mathbb{R}_-$ can be similarly dealt with.
Applying (\ref{equ416}) on $I_n^1$, it follows that
\begin{align*}&\sup_{t\in I_n^1}\|r_n^l\|_{L^2}+\|r_n^l\|_{L_{t,x}^3[I_n^1]}
\lesssim \epsilon\|r_n^l\|_{L_{t,x}^3[I_n^1]}
+\|r_n^l\|_{L_{t,x}^3[I_n^1]}^3+\|r_n^l\|_{L_{t,x}^3[I_n^1]}^2+\|r_n^l(0,\cdot)\|_{L^2}
+o_n(1).\end{align*} By choosing $\epsilon$ sufficiently small, we
obtain
$$\sup_{t\in I_n^1}\|r_n^l\|_{L^2}+\|r_n^l\|_{L_{t,x}^3[I_n^1]}
\lesssim
\|r_n^l(0,\cdot)\|_{L^2}+\sum_{\alpha=2}^3\|r_n^l\|_{L_{t,x}^3[I_n^1]}^\alpha+o(1).$$
Observe that, by the definition of the nonlinear profile $U_n^j$, we
have $$\lim_{n\rightarrow\infty}\|r_n^l(0,\cdot)\|_{L^2}=0$$ for
every $l\geq1$. This fact and a standard bootstrap argument show
easily that $$\lim_{n\rightarrow\infty}\Big(\sup_{t\in
I_n^1}\|r_n^l\|_{L^2}+\|r_n^l\|_{L_{t,x}^3[I_n^1]}\Big)\buildrel{l\rightarrow\infty}\over{\longrightarrow}0.$$
This gives, in particular
$$\lim_{n\rightarrow\infty}\|r_n^l(a_n^1,\cdot)\|_{L^2}\buildrel{l\rightarrow\infty}\over{\longrightarrow}0$$
and allows us to repeat the same argument on $I_n^2$. We iterate the
same process for every $1\leq i\leq p$. Since $I=I_n^1\cup
I_n^2\cup\cdots\cup I_n^p$ and $p$ is finite independently of $n$
and $l$, we get
$$\lim_{n\rightarrow\infty}\big(\|r_n^l\|_{L_{t,x}^3[I_n]}+\sup_{t\in
I_n}\|r_n^l\|_{L^2}\big)\rightarrow0$$ as $l\rightarrow\infty$, which is (\ref{equ64}). \\

{\bf Step 2:} Now we prove the equivalence of (i) and (ii).

$(i)\Rightarrow(ii)$:

Suppose that for all $j$,
$\displaystyle\lim_{n\rightarrow\infty}\|\Gamma_n^jU^j\|_{L_{t,x}^3[I_n]}<+\infty$,
then
$$\big\|u_n\big\|_{L_{t,x}^3[I_n]}\leq\sum_{j=1}^l\big\|U_n^j\big\|_{L_{t,x}^3[I_n]}
+\big\|r_n^l\big\|_{L_{t,x}^3[I_n]}+\big\|w_n^l\big\|_{L_{t,x}^3[I_n]}.$$
From (\ref{p2}), we have
$$\limsup_{n\rightarrow\infty}\|w_n^l\|_{L_{t,x}^3[I_n]}\xrightarrow{l\rightarrow\infty}0
\quad\text{and}\quad\lim_{n\rightarrow\infty}\|r_n^l\|_{L_{t,x}^3[I_n]}\xrightarrow{l\rightarrow\infty}0
.$$ It immediately follows that
$$\lim_{n\rightarrow\infty}\|u_n\|_{L_{t,x}^3[I_n]}<+\infty.$$

$(ii)\Rightarrow(i)$:

 If (i) does not hold, there
exists a family of $\tilde{I}_n\subset I_n$ with 0 included, such
that
$$\sum_{j=1}^\infty\lim_{n\rightarrow\infty}\big\|U_n^j\big\|_{L_{t,x}^3[\tilde{I}_n]}^3>M$$
for arbitrary large $M$ and
$$\|u_n\|_{L_{t,x}^3[\tilde{I}_n]}<\infty.$$

By the orthogonality, we have
$$\lim_{n\rightarrow\infty}\|u_n\|_{L_{t,x}^3[\tilde{I}_n]}^3
\geq\sum_{j=1}^\infty\lim_{n\rightarrow\infty}\|U_n^j\|^3_{L_{t,x}^3[\tilde{I}_n]}>M.$$
This leads to
$$\lim_{n\rightarrow\infty}\|u_n\|_{L_{t,x}^3[I_n]}^3\geq
\lim_{n\rightarrow\infty}\|u_n\|^3_{L_{t,x}^3[\tilde{I}_n]}>M,$$
which implies that
$$\lim_{n\rightarrow\infty}\|u_n\|_{L_{t,x}^3[I_n]}=+\infty.$$
This contradicts (ii). This completes the proof of Theorem
\ref{thm41}.

\vskip0.2cm
 {\noindent \it Proof of Theorem \ref{thm-mass1}.} We choose
$\{u_{0,n}\}$ such that $\|u_{0,n}\|_{L^2}\downarrow \delta_0$, let
$u_n$ is the solution of (\ref{equ42}) with data $u_{0,n}$. By the
definition of $\delta_0$, we can assume that the interval of
existence for $u_n$ is finite. By time translation and scaling, we
may assume that $\{u_n\}_{n=1}^\infty$ is well defined on $[0,1]$,
and
$$\lim_{n\rightarrow\infty}\|u_n\|_{L_t^3([0,1],L_x^3)}=+\infty.$$
Let $\{U^j, V^j, \rho_n^j, s_n^j, \xi_n^j, x_n^j\}$ be the family of
linear and nonlinear profiles associated to $\{u_n\}_{n=1}^\infty$
via Lemma \ref{lem44} and Theorem \ref{thm41}. Then the equivalence
in Theorem \ref{thm41} implies that there exists a $j_0$ such that
$U^{j_0}$ blows up. On one hand, by the definition of
$B_{\delta_0}$,
$$\|V^{j_0}\|_{L^2}\geq\delta_0.$$ On the other hand, we have
$$\sum_{j\geq0}\|V^{j_0}\|_{L^2}^2\leq\lim_{n\rightarrow\infty}\|u_{0,n}\|_{L^2}^2=\delta_0^2.$$
Thus by mass conservation and the definition of nonlinear profile,
we have
$$\|U^{j_0}\|_{L^2}=\|V^{j_0}\|_{L^2}\leq\delta_0.$$ Therefore,
$$\|U^{j_0}\|_{L^2}=\delta_0.$$ Because $U^{j_0}$ is the solution of
(\ref{equ42}) satisfying $U(s^{j_0},x)=V(s^{j_0},x)$, where
$s^{j_0}=\lim_{n\rightarrow\infty} s^{j_0}_n$. If $s^{j_0}$ is
finite, then $U^{j_0}$ is the blow up solution with minimal mass. If
$s^{j_0}=\infty$, we can use the pseudo-conformal transformation to
get a blow up solution with minimal mass. This shows the existence
of initial data such that solution of (\ref{equ42}) blows up in
finite time for $t>0$. In the proof of Theorem \ref{thm-mass2} we
will show that there exists an initial data $u_0\in
L^2(\mathbb{R}^4)$ with $\|u_0\|_{L^2}=\delta_0$, such that the
solution u of (\ref{equ42}) blows up for both $t>0$ and $t<0$.

\vskip0.2cm
 {\noindent \it Proof of Theorem \ref{thm-mass2}}.  (i) Suppose $u$
is a solutions of (\ref{equ42}) which blows up at finite time
$T^*>0$ and $\{t_n\}_{n=1}^\infty$ is a sequence increasingly going
to $T^*$ as $n\rightarrow\infty$. Let
$$u_n(t,x)=u(t_n+t,x),$$
then $\{u_n\}_{n=1}^\infty$ is a family of
solutions on $I_n=[-t_n,T^*-t_n)$. Moreover, we have
$$\lim_{n\rightarrow\infty}\Big\|u_n\Big\|_{L_{t,x}^3\in[0,T^*-t_n)}=
\lim_{n\rightarrow\infty}\Big\|u_n\Big\|_{L_{t,x}^3\in[-t_n,0]}=+\infty.$$
Since $\|u_n\|_{L^2}$ is bounded due to $L^2$ conservation, we can
apply Lemma \ref{lem44} and then Theorem \ref{thm41} on
$I_n=[0,T^*-t_n)$ to get that there exists some $j_0$ such that the
nonlinear profile $\{U^{j_0}, \rho_n^{j_0},s_n^{j_0},\xi_n^{j_0},
x_n^{j_0}\}$ satisfies
\begin{equation}\label{equ417}\lim_{n\rightarrow\infty}\big\|U^{j_0}\big\|_{L_{t,x}^3[I_n^{j_0}]}=+\infty,\end{equation}
where $$I_n^{j_0}:=[s_n^{j_0},
(\rho_n^{j_0})^2(T^*-t_n)+s_n^{j_0}).$$ In fact, let
$s^{j_0}=\lim_{n\rightarrow\infty}s_n^{j_0}$, then
$s^{j_0}\neq\infty$, otherwise, $I_{n}^{j_0}\rightarrow\emptyset$
and (\ref{equ417}) is impossible. This implies either
$s^{j_0}=-\infty$ or $s^{j_0}=0$ (up to translation). If
$s^{j_0}=0$, let $U^{j_0}$ be the solution of (\ref{equ41}) with
initial data $V^{j_0}$, then (\ref{equ417}) implies $U^{j_0}$ blows
up at time $T^*_{j_0}\in(0,+\infty)$ and
\begin{equation}\label{equ418}\lim_{n\rightarrow\infty}(\rho_n^{j_0})^2(T^*-t_n)\geq
T_{j_0}^*.\end{equation} If we assume also that
$\|u_0\|_{L^2}<\sqrt{2}\delta_0$, then there is at most one linear
profile with $L^2$-norm greater than $\delta_0$ thanks to
(\ref{p3}). That means that the profile $U^{j_0}$ founded above is
the only blow up nonlinear profile (since all the other profiles
have $L^2$ norm less than $\delta_0$ and then they are global ).
By repeating the same argument in $I_n=[-t_n, 0]$, we get
$$\lim_{n\rightarrow\infty}\|U^{j_0}\|_{L_{t,x}^3[I_n^{j_0}]}=+\infty,
\ \ I_n^{j_0}=[-(\rho_n^{j_0})^2t_n+s_n^{j_0}, s_n^{j_0}].$$ This
implies that $s^{j_0}\neq-\infty$. Hence $s^{j_0}=0$ and the
solution $U^{j_0}$ of (\ref{equ42}) with initial data
$V^{j_0}(0,\cdot)$ blows up also for $t<0$. Thus the nonlinear
profile $U^{j_0}$ is the solution of (\ref{equ42}) which blows up
for both $t<0$ and $t>0$.

(ii) The linear decomposition yields
$$(\mathbf{\Gamma}_n^{j_0})^{-1}(e^{it\Delta}(u(t_n,\cdot))=V^{j_0}+
\sum_{1\leq j\leq l; j\neq
j_0}(\mathbf{\Gamma}_n^{j_0})^{-1}\mathbf{\Gamma}_n^jV^j+(\mathbf{\Gamma}_n^{j_0})^{-1}w_n^l.$$
The family $\{\mathbf{\Gamma}_n^j\}_{j=1}^\infty$ is pairwise
orthogonal, so for every $j\neq j_0$,
$$(\mathbf{\Gamma}_n^{j_0})^{-1}\mathbf{\Gamma}_n^{j}V^j\Extend{-}{-}{\rightharpoonup}{}{n\rightarrow\infty}0
\ \text{weakly in}\ L^2.$$ Then
$$(\mathbf{\Gamma}_n^{j_0})^{-1}(e^{it\Delta}(u(t_n,\cdot))\Extend{-}{-}{\rightharpoonup}{}{n\rightarrow\infty} V^{j_0}+\tilde{w}^l\ \ \text{weakly}.$$
where $\tilde{w}^l$ denote the weak limit of
$(\mathbf{\Gamma}_n^{j_0})^{-1}w_n^l$. However, we have
$$\|\tilde{w}^l\|_{L_{t,x}^3}\leq\lim_{n\rightarrow\infty}\|w^l_n\|_{L_{t,x}^3}\xrightarrow{l\rightarrow+\infty}0.$$
By the uniqueness of weak limit, we get $\tilde{w}^l=0$ for every
$l\geq j_0$. Hence, we obtain
$$(\mathbf{\Gamma}_n^{j_0})^{-1}(e^{it\Delta}(u(t_n,\cdot))\Extend{-}{-}{\rightharpoonup}{}{n\rightarrow\infty}
V^{j_0}.$$

We need the following lemma:
\begin{lemma}[\cite{MeV}]\label{lem46}Let $\{\varphi_n\}_{n\geq 1}$ and $\varphi$
be in $L^2(\mathbb{R}^4)$. The following statement is equivalent:

(1) $\varphi_n\rightharpoonup\varphi$ weakly in $L^2(\mathbb{R}^4)$.

(2) $e^{it\Delta}\varphi_n\rightharpoonup e^{it\Delta}\varphi$ in
$L_{t,x}^3(\mathbb{R}^{4+1})$\end{lemma}

Applying this lemma to
$(\mathbf{\Gamma}_n^{j_0})^{-1}(e^{it\Delta}(u(t_n,\cdot))$, we get
$$e^{-is_n\Delta}\Big(\rho_n^2e^{ix\cdot\xi_n}e^{i\theta_n}u(t_n,\rho_nx+x_n)\Big)\rightharpoonup V^{j_0}(0,\cdot)$$
with $$s_n=s_n^{j_0},\ \ \rho_n=\frac{1}{\rho_n^{j_0}},\ \
\theta_n=\frac{x_n^{j_0}\xi_n^{j_0}}{\rho_n^{j_0}},\ \
x_n=\frac{-x_n^{j_0}}{\rho_n^{j_0}},\ \
\xi_n=-\frac{\xi_n^{j_0}}{\rho_n^{j_0}}.$$ Up to subsequence, we
can assume that $e^{i\theta_n}\rightarrow e^{i\theta}$. Since
$s_n\rightarrow0$, we get
\begin{equation}\label{equ65}
\rho_n^2e^{ix\cdot\xi_n}u(t_n,\rho_nx+x_n)\rightharpoonup
e^{-i\theta}V^{j_0}(0,\cdot).
\end{equation}
The associated
solution is $e^{-i\theta}U^{j_0}$. (\ref{equ418}) gives
$$\lim_{n\rightarrow\infty}\frac{\rho_n}{\sqrt{T^*-t_n}}\leq\frac{1}{\sqrt{T_{j_0}^*}}.$$
This completes the proof of Theorem \ref{thm-mass2}.

(iii) Let $u$ be a solution of (\ref{equ1}) with
$\|u_0\|_{L^2}<\sqrt{2}\delta_0$, which blows up at finite time
$T^*>0$. Let $\{t_n\}_{n=1}^\infty$ be any time sequence such that
 $t_n\uparrow T^*$ as $n\rightarrow\infty$. So there exist $V\in
L^2(\mathbb{R}^4)$ with $\|V\|_{L^2}\geq\delta_0$ and a sequence
$\{\rho_n,\xi_n,x_n\}\subset\mathbb{R}_+^*\times\mathbb{R}^4\times\mathbb{R}^4$
such that up to a subsequence,
$$(\rho_n)^2e^{ix\cdot\xi_n}u(t_n,\rho_nx+x_n)\Extend{-}{-}{\rightharpoonup}{}{n\rightarrow\infty}V$$
and $$\lim_{n\rightarrow\infty}\frac{\rho_n}{\sqrt{T^*-t_n}}\leq A$$
for some $A\geq0$. Thus we have
$$\lim_{n\rightarrow\infty}\rho_n^4\int_{|x|\leq
R}|u(t_n,\rho_nx+x_n)|^2dx\geq\int_{|x|\leq R}|V|^2dx$$ for every
$R\geq0$.  This implies that
$$\lim_{n\rightarrow\infty}\sup_{y\in\mathbb{R}^4}\int_{|x-y|\leq
R\rho_n}|u(t_n,x)|^2dx\geq\int_{|x|\leq R}|V|^2dx.$$ Since
$\frac{\sqrt{T^*-t}}{\lambda(t)}\rightarrow0$ as $t\uparrow T^*$,
it follows that $\frac{\rho_n}{\lambda(t_n)}\rightarrow0$ and then
$$\lim_{n\rightarrow\infty}\sup_{y\in\mathbb{R}^4}\int_{|x-y|\leq
\lambda(t_n)}|u(t_n,x)|^2dx\geq\int|V|^2dx\geq\delta_0^2.$$
 Since $\{t_n\}_{n=1}^\infty$ is an arbitrary sequence, we
infer
$$\liminf_{t\rightarrow T}\sup_{y\in\mathbb{R}^4}\int_{|x-y|\leq
\lambda(t)}|u(t,x)|^2dx\geq\delta_0^2.$$ However for every
$t\in[0,T)$, the function
$y\mapsto\int_{|x-y|\leq\lambda(t)}|u(t,x)|^2dx$ is continuous and
goes to 0 at infinity. As a consequence, we get
$$\sup_{y\in\mathbb{R}^4}\int_{|x-y|\leq\lambda(t)}|u(t,x)|^2dx=
\int_{|x-x(t)|\leq\lambda(t)}|u(t,x)|^2dx$$ for some
$x(t)\in\mathbb{R}^4$ and this completes the proof of Theorem
\ref{thm-mass2}.

\vskip 0.2cm
 {\noindent \it Proof of Corrolary \ref{co}.} In context of the
proof of Theorem \ref{thm-mass2} we assume also that
$$\|u_n\|_{L^2}=\|u_0\|_{L^2}=\delta_0.$$
(\ref{p3}) gives that $$\|V^{j_0}\|_{L^2}\leq\delta_0.$$ It
follows that $$\|V^{j_0}\|_{L^2}=\delta_0.$$ This implies that
there exists a unique profile $V^{j_0}$ and the weak limit in
(\ref{equ65}) is strong.

\vskip 0.2cm
 {\bf Acknowledgement: }
The authors  thank professor S. Keraani for several lectures on
his work \cite{HMKer} and helpful discussions, and thank Pin Yu
for helpful discussion.  C. Miao and G. Xu were partly supported
by the NSF of China (No.10725102, No.10726053), and L. Zhao was
supported by China postdoctoral science foundation project.
\vskip0.3cm


\begin{center}

\end{center}

\begin{thebibliography}{99}
\addcontentsline{toc}{section}{References}
\bibitem{Ba}V. Banica, Remarks on the blow-up for the Schr\"{o}dinger
equation with critical mass on a plane domain. Ann. Sc. Norm. Super.
Pisa Cl. Sci. (5) 3 (2004), no. 1, 139-170.


\bibitem{BeV} P. B\'{e}gout and A. Vargas, Mass concentration phenomena for the
$L^2$-critical nonlinear Schr\"{o}dinger equaiton. Trans. Amer.
Math. Soc., 359(2007), 5257-5282.

\bibitem{Bo}J. Bourgain, Refinements of Strichartz inequaltiy and
applications to 2D-NLS with critical nonlinearity. IMRN, 8(1998)
253-283.

\bibitem{CaK}R. Carles and S. Keraani, On the role of quadratic
oscillations in nonlinear Schr\"{o}dinger equation II. The
$L^2$-critical case. Trans. Amer. Math. Soc., 359(2007), 33-62.

\bibitem{Ca03}T. Cazenave, {\it Semilinear Schr\"{o}dinger equations}. Courant
Lecture Notes in Mathematics, Vol. 10. New York: New York University
Courant Institute of Mathematical Sciences, 2003.

\bibitem{FrL04}J. Fr\"{o}hlich and E. Lenzmann, Mean-field limit of quantum
Bose gases and nonlinear Hartree equation. S\'{e}minaire:
\'{E}quations aux D\'{e}riv\'{e}es Partielles. 2003-2004, S\'{e}min.
\'{E}qu. D\'{e}riv. Partielles, \'{E}cole Polytech., Palaiseau,
2004, pp. Exp. No. XIX, 26.

\bibitem{Gi}J. Ginibre, {\it Introduction aux \'{e}quations de
Schr\"odinger non lin\'eaires.} Master course, 94-95.

\bibitem{GiV}J. Ginibre and G. Velo, Scattering theory in the
energy space for a class of Hartree equations. Nonlinear wave
equations (Providence, RI, 1998), 29-60, Contemp. Math., 263, Amer.
Math. Soc., Providence, RI, 2000.


\bibitem{HMKer}T. Hmidi and S. Keraani, Blowup theory for the
critical nonlinear Schr\"{o}dinger equations revisited. IMRN,
46(2005) , 2815-2828.

\bibitem{KeT98}M. Keel and T. Tao, Endpoint Strichartz estimates. Amer. J.
Math., 120:5(1998), 955-980.


\bibitem{Ker1}S. Keraani, On the defect of compactness for the
Strichartz estimates of the Schr\"{o}dinger equations. J. Diff.
Equa., 175(2001) 353-392.

\bibitem{Ker2}S. Keraani, On the blow up phenomenon of the
critical nonlinear Schr\"{o}dinger equation. J. Funct. Anal.,
235(2006), 171-192.

\bibitem{KrLR08}J. Krieger, E. Lenzmann and P. Raphael, On stability
of pseudo-conformal blowup for $L^2$-critical Hartree equation.
arXiv:0808.2324.

\bibitem{kong}M. K. Kwong, Uniequeness of positive solutions of $\Delta u-u+u^p=0$ in
$\Bbb R^n$. Arch. Rat. Mech. Anal., 105(1989), 243-266.


\bibitem{LiMZ08}D. Li, C. Miao and X. Zhang, The focusing energy-critical Hartree equation. To appear in J. Diff.
Equa..

\bibitem{lieb}E. H. Lieb, Existence and uniqueness of the minimizing
solution of Choquar's nonlinear equation. Stud. Appl. Math., 57
(1977), 93-105.

%

\bibitem{Liu}S. Liu, Regularity, symmetry, and uniqueness of some
integral type quasilinear equations with nonlocal nonlinearities.
Preprint.

\bibitem{Mer1}F. Merle, Blow-up phenomena for critical nonlinear
Schr\"{o}dinger and Zakharov equations. Proceeding of the
International Congress of Mathematicians (Berlin, 1998), Doc. Math.
extra. Vol. III(1998), 57-66.

\bibitem{Mer2}F. Merle, Determination of blow-up solutions with
minimal mass for nonlinear Schr\"{o}dinger equaitons with critical
power. Duke Math. J., 69:2(1993), 427-454.

\bibitem{MeT}F. Merle and Y. Tsutsumi, $L^2$ concentration of
blow-up solutions for the nonlinear Schr\"{o}dinger equation with
critical power nonlinearity. J. Diff. Equa., 84(1990), 205-214.

\bibitem{MeV}F. Merle and L. Vega, Compactness at blow-up time for
$L^2$ solutions of the critical nonlinear Schr\"{o}dinger equation
in 2D. IMRN, 8(1998), 399-425.

\bibitem{MeR1}F. Merle and P. Raphael, Sharp upper bound on the
blow-up rate for the critical nonlinear Schr\"{o}dinger equation.
GAFA, 13(2003), 591-642.

\bibitem{MeR2}F. Merle and P. Raphael, On universality of blow-up
profile for $L^2$ critical nonlinear Schr\"{o}dinger equation.
Invent. Math., 156 (2004), 565-672.

\bibitem{MeR05}F. Merle and P. Raphael, On a sharp lower bound on
the blow-up rate for the $L^2$ critical nonlinear Schr\"{o}dinger
equation. J. Amer. Math. Soc., 19:1(2005), 37-90.

\bibitem{MiXZ1}C. Miao, G. Xu and L. Zhao, The Cauchy problem of the Hartree
equation. J. PDE, 21(2008), 22-44.

\bibitem{MiXZ2}C. Miao, G. Xu and L. Zhao, Global well-posedness and scattering
 for the energy-critical, defocusing Hartree equation for radial
 data. J. Funct. Anal., 253(2007), 605-627.

\bibitem{MiXZ3}C. Miao, G. Xu and L. Zhao, Global well-posedness and scattering
for the energy-critical, defocusing Hartree equation in
$\mathbb{R}^{1+n}$. Preprint.

\bibitem{MiXZ5}C. Miao, G. Xu and L. Zhao, Global well-posedness,
scattering and blow-up for the energy-critical, focusing Hartree
equation in the radial case. To appear in Coll. Math.

\bibitem{MiXZ4}C. Miao, G. Xu and L. Zhao, Global well-posedness and
scattering for the mass-critical Hartree equation with radial data.
To appear in J. Math. Pures Appl..

\bibitem{Na}K. Nakanishi, Energy scattering for Hartree equations.
Math. Res. Lett., 6(1999), 107-118.

 \bibitem{nawa}H. Nawa, ``Mass concentration'' phenomenon for the
 nonlinear Schr\"{o}dinger equation with the critical power
 nonlinearity. Fukcial. Ekvac., 35:1(1992), 1-18.

\bibitem{tao1}T. Tao, M. Visan, and X. Zhang, Minimal-mass blowup solutions of the mass-critical
NLS. To appear in Forum Math.

\bibitem{tao2}R. Killip, T. Tao and M. Visan, The cubic nonlinear
Schr\"odinger equation in two dimensions with radial data. Preprint.

\bibitem{tao3} R. Killip, M. Visan and X. Zhang, The mass-critical
nonlinear Schr\"odinger equation with radial data in dimensions
three and higher. Preprint.


\bibitem{wein89}M. Weinstein, The Nonlinear Schr\"{o}dinger
Equation-Singularity Formation, Stability and Dispersion. pp.
213-232 in: The connection between infinite-dimensional and
finite-dimensional dynamical systems, Contemporary Math., No. 99,
Amer. Math. Soc., Providence, R. I., 1989.

\bibitem{wiki}http://tosio.math.toronto.edu/wiki/index.php/Hartree equation.
\end{thebibliography}
\end{document}